\documentclass{article}

\usepackage[a4paper,top=4.0cm,bottom=2.0cm,left=2.0cm,right=2.0cm]{geometry}

\usepackage{amsmath,amsfonts,mathrsfs,amsfonts}
\usepackage{epsfig,graphicx}
\usepackage{color}

\newtheorem{lemma}{Lemma}[section]
\newtheorem{definition}{Definition}[section]
\newtheorem{remark}{Remark}[section]

\begin{document}

\title{Second-order accurate ensemble transform particle filters}

\author{Walter Acevedo \thanks{Universit\"at Potsdam, 
Institut f\"ur Mathematik, Karl-Liebknecht-Str. 24/25, D-14476 Potsdam, Germany}  \and
Jana de Wiljes\thanks{Universit\"at Potsdam, 
Institut f\"ur Mathematik, Karl-Liebknecht-Str. 24/25, D-14476 Potsdam, Germany} 
\and Sebastian Reich\thanks{Universit\"at Potsdam, 
Institut f\"ur Mathematik, Karl-Liebknecht-Str. 24/25, D-14476 Potsdam, Germany ({\tt sreich@math.uni-potsdam.de}) and University of Reading, Department of Mathematics and Statistics, Whiteknights, PO Box 220, Reading RG6 6AX, UK.} 
}

\maketitle

\begin{abstract} Particle filters (also called sequential Monte Carlo methods) are widely used for
state and parameter estimation problems in the context of nonlinear evolution equations. 
The recently proposed ensemble transform particle filter (ETPF) 
(S.~Reich, {\it A non-parametric ensemble transform
method for Bayesian inference}, SIAM J.~Sci.~Comput., 35, (2013), pp. A2013--A2014) replaces the
resampling step of a standard particle filter by a linear transformation which allows for a hybridization of
particle filters with ensemble Kalman filters and renders the resulting hybrid filters applicable to spatially 
extended systems. However, the linear transformation step is computationally expensive
and leads to an underestimation of the ensemble spread for small and moderate ensemble sizes.
Here we address both of these shortcomings by developing second-order accurate extensions of the
ETPF. These extensions allow one in particular to replace the exact solution of a linear transport problem by its Sinkhorn approximation. It is also demonstrated that the nonlinear ensemble 
transform filter (NETF) arises as a special case of our general framework. We illustrate
the performance of the second-order accurate filters for the chaotic Lorenz-63 and 
Lorenz-96 models and a dynamic scene-viewing model. The numerical results for the Lorenz-63 and
Lorenz-96 models demonstrate that significant accuracy improvements can be achieved in comparison
to a standard ensemble Kalman filter and the ETPF for small to moderate ensemble sizes. 
The numerical results for the scene-viewing model reveal, on the other hand, that second-order
corrections can lead to statistically inconsistent samples from the posterior parameter distribution.
\end{abstract}

\noindent {\bf AMS(MOS) subject classifications.} 65C05, 62M20, 93E11, 62F15, 86A22

%
\section{Introduction}
%

Data assimilation (DA) denotes the broad topic of combining evolution models with partial observations
of the underlying dynamical process \cite{sr:evensen,sr:stuart15,sr:reichcotter15}. 
DA algorithms come in the form of variational and/or ensemble-based
methods \cite{sr:stuart15}. In this paper, we focus on ensemble-based DA methods and their robust and
efficient implementation. The ensemble Kalman filter (EnKF) \cite{sr:evensen} is by far the most popular 
ensemble-based DA method and has found widespread application in the geosciences. 
However, EnKFs lead to inconsistent approximations for partially observed nonlinear processes. On the contrary,
particle filters (PF) (also called sequential Monte Carlo methods) \cite{sr:Doucet} 
lead to consistent approximations but  typically require ensemble sizes much larger than those required 
for EnKFs in order to track the underlying reference process \cite{sr:bengtsson08}. 

In order to overcome these shortcomings, we are currently witnessing a strong trend towards hybrid filters 
which combine EnKFs with PFs and which are applicable to strongly nonlinear systems 
under small or moderate ensemble sizes. We mention here the Gaussian mixture filters 
(such as, for example, \cite{sr:stordal11}), the rank histogram filter \cite{sr:anderson10,sr:snyder13}, moment matching ensemble filters \cite{sr:xiong06,sr:lei11,sr:toedter15}, the ensemble Kalman particle filter \cite{sr:frei13}, and the hybrid ensemble transform particle filter \cite{sr:CRR15}.

In this paper, we focus on improved implementations of the ensemble transform particle filter (ETPF)
\cite{sr:reich13,sr:reichcotter15} and its hybridrization with the EnKF \cite{sr:CRR15}. The ETPF requires the solution of a linear transport problem in each assimilation step, which 
renders the methods substantially
more expensive than an EnKF. Computationally attractive alternatives, such as the Sinkhorn approximation \cite{sr:cuturi13}, lead to unstable implementations since the ensemble becomes underdispersive. We address
this problem by introducing a variant of the ETPF, which is second-order accurate independent of the actual
solution procedure for the underlying optimal transport problem. An ensemble filter is called second-order accurate
if the posterior mean and covariance matrix of the ensemble are in agreement with their importance sampling 
estimates from a Bayesian inference step. Second-order accurate particle filters have first been proposed in \cite{sr:xiong06} and since 
then several variants of it have been developed \cite{sr:lei11,sr:toedter15}. Here we instead consider second-order corrections to the ETPF. Such corrections
require the solution of a continuous-time algebraic Riccati equation 
\cite{sr:Wonham68,sr:laub79}. 
The correction term vanishes as the ensemble size approaches infinity in agreement with the consistency of the ETPF \cite{sr:reich13}. 

The paper is organized as follows. The general framework of ensemble transform filters is 
summarized in Section \ref{sec:FORDA}. Section \ref{sec:2ndETPF} summarizes the ETPF and
introduces the second-order correction step. Numerical solution procedures
for the associated continuous-time algebraic Riccati equation are discussed in 
Section{sec:Riccati}. The Sinkhorn approximation to
the optimal transport problem of the ETPF is introduced in Section \ref{sec:sinkhorn} and the overall
second-order accurate implementation of the ETPF is summarized in Section \ref{sec:AS}. 
Numerical results are provided in Section \ref{sec:NE}, where the behavior of the new method is
demonstrated  for the highly nonlinear and chaotic
Lorenz-63 \cite{sr:lorenz63} and Lorenz-96 \cite{sr:lorenz96} 
models. Here we repeat the experiments from \cite{sr:CRR15} with
the ETPF being replaced by a second-order accurate
variant based on the Sinkhorn approximation to the underlying optimal transport problem. 
We finally also demonstrate the behavior of the new filters for parameter estimation of
the scene-viewing model {\it SceneWalk} \cite{sr:Engbert2015}.

%
\section{Ensemble-based forecasting-data assimilation systems} \label{sec:FORDA}
%

Let us assume that observations ${\bf y}^{\rm obs}(t_k)\in \mathbb{R}^{N_y}$ 
become available at time instances 
$t_k$, $k=1,\ldots,K$, and are related to the state variables ${\bf z} \in \mathbb{R}^{N_z}$ 
of an evolution model 
\begin{equation} \label{model}
{\bf z}(t_k) = {\cal M}({\bf z}(t_{k-1}))
\end{equation}
via the likelihood function
\begin{equation} \label{likelihood}
\pi({\bf y}|{\bf z}) = \frac{1}{(2\pi)^{N_y/2} |{\bf R}|^{1/2}} \exp \left( -\frac{1}{2} (h({\bf z})-{\bf y})^{\rm T} 
{\bf R}^{-1} (h({\bf z})-{\bf y})\right),
\end{equation}
where ${\bf R}\in \mathbb{R}^{N_y\times N_y}$ denotes the measurement error covariance matrix.

An ensemble-based forecasting-data assimilation (FOR-DA) systems will produce two sets of ensembles
of size $M$ at any $t_k$. First we have the forecast ensemble $\{{\bf z}_i^{\rm f}\}_{i=1}^M$
which approximates the conditional distribution $\pi({\bf z},t_k|{\bf y}^{\rm obs}_{1:k-1})$ and, second,
we have the analysis ensemble $\{{\bf z}_i^{\rm a}\}_{i=1}^M$, which approximates the conditional
distribution $\pi({\bf z},t_k | {\bf y}^{\rm obs}_{1:k})$. Here 
\begin{equation}
{\bf y}^{\rm obs}_{1:l} = ({\bf y}^{\rm obs}(t_1),{\bf y}^{\rm obs}(t_2),\ldots,{\bf y}^{\rm obs}(t_l)) \in \mathbb{R}^{N_y \times l}
\end{equation}
denotes the complete set of observations from $t=t_1$ to $t=t_l$. Also note that
\begin{equation} \label{model_step}
{\bf z}_i^{\rm f}(t_k) = {\cal M}({\bf z}_i^{\rm a}(t_{k-1}))
\end{equation}
and that FOR-DA systems primarily differ  in the employed data assimilation
algorithms. 

The data assimilation algorithms considered in this paper are all of the form of 
a linear ensemble transform filter (LETF) \cite{sr:reichcotter15}:
\begin{equation} \label{transform}
{\bf z}_j^{\rm a}(t_k) = \sum_{i=1}^M {\bf z}_i^{\rm f}(t_k) \,d_{ij}(t_k)
\end{equation}
where the entries $d_{ij}(t_k)$ of the $M\times M$ transformation matrix ${\bf D}(t_k) = \{d_{ij}(t_k)\}$ 
are subject to the constraint
\begin{equation} \label{constrainttransformationmatrixentries}
\sum_{i=1}^M d_{ij}(t_k) = 1
\end{equation}
for all $j=1,\ldots,M$.
In other words, provided that $M\le N_z$, the members ${\bf z}_j^{\rm a}(t_k)$ of the analysis ensemble lie in the $(M-1)$-dimensional hyperplane spanned by the forecast ensemble ${\bf z}_i^{\rm f}(t_k)$ with $i\in\{1,\dots,M\}$. Note that the entries of ${\bf D}$ can be negative. See \cite{sr:reichcotter15}.
One well known exemplary class of DA algorithms that has the LETF structure is the family of EnKFs \cite{sr:evensen,sr:reichcotter15}. It has long been acknowledged  that EnKFs are very robust yet the underlying Gaussianity and linearity assumptions limit their applicability to more general systems. To address the shortcomings of traditional techniques such as the EnKFs, other algorithms that are applicable for nonlinear model scenarios and are computational feasible when employed to high-dimensional systems have been proposed. For example, the nonlinear ensemble transform filter (NETF) of \cite{sr:xiong06,sr:toedter15} provides an example of a 
particle filter in the form of an LETF  based upon the normalized importance weights
\begin{equation} \label{nweights}
w_i(t_k) := \frac{\widehat{w}_i(t_k)}{\sum_{j=1}^M \widehat{w}_j(t_k)}
\end{equation}
with $\widehat{w}_i(t_k) =  \pi({\bf y}(t_k)|{\bf z}_i^{\rm f}(t_k))$, which reproduces the first and second-order moments of the posterior distribution. 
The main focus of this paper is, however, on the ETPF which can also be formulated in the form of (\ref{transform}) \cite{sr:evensen,sr:reichcotter15} with ${\bf D}(t_k) = \{d_{ij}(t_k)\}$ being defined via minimization of the cost functional 
\begin{equation} \label{cost}
J({\bf D}(t_k)) = \sum_{i,j=1}^M d_{ij}(t_k) \|{\bf z}_i^{\rm f}(t_k) - {\bf z}_j^{\rm f}(t_k)\|^2
\end{equation}
subject to $d_{ij}(t_k) \ge 0$,  (\ref{constrainttransformationmatrixentries}) and
\begin{equation}\label{constrainttransformationmatrixentries2}
\frac{1}{M} \sum_{j=1}^M d_{ij}(t_k) = w_i(t_k)
\end{equation}
\cite{sr:reich13}. The key idea of the ETPF is to approximate a transfer map between the random variable, $Z^{\rm f}(t_k)$, distributed according to $
\pi_{Z^{\rm f}}({\bf z},t_k)$ and the random variable, $Z^{\rm a}(t_k)$, associated with $\pi_{Z^{\rm a}}({\bf z},t_k)$. 
This map induces a coupling of the respective densities that is optimal in the sense that it minimizes the expected distance between the two random variables, i.e.,
\begin{equation}\label{Monge-Kantorovitch problem}
\mu^*_Z=\arg \inf_{\mu\in\Pi(\pi_{Z^{\rm f}},\pi_{Z^{\rm a}})}\sqrt{\mathbb{E}||Z^{\rm f}(t_k)-
Z^{\rm a}(t_k)||^2}.
\end{equation}
Intuitively it is clear that the correlation between the forecast and the analysis random variable is increased via optimization of (\ref{Monge-Kantorovitch problem}) 
and thus creates a strong relation between the prior and the posterior. Since we only rely on importance weights, our filter is also applicable to non-Gaussian likelihood functions. 
The ETPF can also be applied to spatially extended systems using the idea of localization \cite{sr:reich15} and has been combined with EnKFs in an hybridization approach \cite{sr:CRR15}. 
While the ETPF convergence to the true posterior distribution in the limit of $M\to \infty$ \cite{sr:reich13}, this is not the case for the EnKF or the NETF, in
general.  However, the ETPF is computationally expensive and underestimates the ensemble spread (covariance matrix) for finite ensemble sizes (see example 8.11 in \cite{sr:reichcotter15}). 
Both of these shortcomings will be addressed by the LETFs proposed in Sections \ref{sec:2ndETPF} and \ref{sec:sinkhorn}.

%
\section{Second-order accurate LETFs} \label{sec:2ndETPF}
%

We now derive second-order accurate LETFs. Here second-order accuracy refers to reproducing the first and second-order moments
exactly according to the importance sampling approach.
\begin{definition}
An LETF (\ref{transform}) is called second-order accurate if 
the analysis mean satisfies
\begin{equation} \label{mean}
\overline{\bf z}^{\rm a}(t_k) = \frac{1}{M}\sum_{i=1}^M {\bf z}^{\rm a}_i(t_k) = 
\sum_{i=1}^M w_i(t_k) {\bf z}_i^{\rm f}(t_k)
\end{equation}
and the analysis covariance matrix
\begin{equation} \label{cma}
\widehat{\bf P}^{\rm a}(t_k) = \frac{1}{M} \sum_{i=1}^M ({\bf z}_i^{\rm a}(t_k)-\overline{\bf z}^{\rm a}(t_k))
({\bf z}_i^{\rm a}(t_k)-\overline{\bf z}^{\rm a}(t_k))^{\rm T} 
\end{equation}
is equal to the covariance matrix defined by the importance weights, i.e.
\begin{equation}\label{cme}
{\bf P}^{\rm a}(t_k) = \sum_{i=1}^M w_i(t_k) ({\bf z}_i^{\rm f}(t_k)-\overline{\bf z}^{\rm a}(t_k))
({\bf z}_i^{\rm f}(t_k)-\overline{\bf z}^{\rm a}(t_k))^{\rm T}.
\end{equation}
\end{definition}

\noindent
\begin{remark}
The covariance matrix (\ref{cme}) derived via importance sampling leads to the denominator $M$ in case of equal weights $w_i = 1/M$. In line with this, the biased version of the empirical covariance (\ref{cma}) is used in this paper. Another option is to introduce the factor 
$\frac{M}{M-1}$ in (\ref{cme}) to obtain the unbiased variant (as is used, for example, in the NETF 
see \cite{sr:toedter15}).
\end{remark}

Since only the DA step of a FOR-DA system is considered in this and the following sections, 
we drop the explicit time-dependence for notational convenience from now on. We introduce the $N_z \times M$ matrix of the forecast ensemble
\begin{equation}
{\bf Z}^{\rm f} = ({\bf z}_1^{\rm f}, {\bf z}^{\rm f}_2,\ldots,{\bf z}_M^{\rm f} ) \in
\mathbb{R}^{N_z\times M}
\end{equation}
and an analog expression
\begin{equation}
{\bf Z}^{\rm a} = ({\bf z}_1^{\rm a}, {\bf z}^{\rm a}_2,\ldots,{\bf z}_M^{\rm a} ) \in
\mathbb{R}^{N_z\times M}
\end{equation}
for the analysis ensemble. Then an LETF (\ref{transform}) can be represented
in the form
\begin{equation} \label{transform_compact}
{\bf Z}^{\rm a} = {\bf Z}^{\rm f} {\bf D}.
\end{equation}
We also introduce the vector ${\bf 1} = (1,1,\ldots,1)^{\rm T} \in \mathbb{R}^{M\times 1}$, 
the vector
\begin{equation}
{\bf w}  = (w_1,\ldots,w_M)^{\rm T} \in \mathbb{R}^{M\times 1}
\end{equation} 
of normalized importance weights (\ref{nweights}), and the diagonal $M \times M$ 
matrix ${\bf W} = \mbox{diag}\,({\bf w})$. Since the analysis mean is provided by (\ref{mean}),
an LETF is first-order accurate if
\begin{equation}\label{LETFfirstorderaccurate}
\frac{1}{M} {\bf Z}^{\rm a} {\bf 1} = {\bf Z}^{\rm f} {\bf w}.
\end{equation}
Equation (\ref{LETFfirstorderaccurate}) holds if ${\bf D}$ satisfies (\ref{constrainttransformationmatrixentries2}), i.e.
\begin{equation}\label{meanconditionequal}
\frac{1}{M} {\bf D}{\bf 1} = {\bf w}.
\end{equation}
Recall that the transformation matrix is also subject to (\ref{constrainttransformationmatrixentries}), which is equivalent to 
${\bf D}^{\rm T} {\bf 1} = {\bf 1}$ \cite{sr:reichcotter15}. In the following, the focus is on first-order accurate LETF characterized by transformation matrices, ${\bf D}$, in the class 
\begin{equation}
\mathcal{D}_1=\{{\bf D} 
\in\mathbb{R}^{M\times M}|\, {\bf D}^{\rm T}{\bf 1} = {\bf 1},~ {\bf D}{\bf 1} = M{\bf w}\, \}.
\end{equation}
These conditions are, for example, satisfied by the transformation matrix 
\begin{equation} \label{simple}
{\bf D}_0 = {\bf w}{\bf 1}^{\rm T},
\end{equation}
which leads to the analysis ensemble 
\begin{equation}
{\bf Z}^{\rm a} = \overline{\bf z}^{\rm a} {\bf 1}^{\rm T}.
\end{equation}

\begin{remark} An EnKF also leads to transformations of the form 
(\ref{transform_compact}) with the associated ${\bf D}_{\rm EnKF}$ satisfying
${\bf D}_{\rm EnKF}^{\rm T} {\bf 1} = {\bf 1}$ but in general not (\ref{meanconditionequal})
\cite{sr:reichcotter15}.
Hence ${\bf D}_{\rm EnKF} \notin \mathcal{D}_1$, in general. The simple modification
\begin{equation}
\widehat{\bf D}_{\rm EnKF} = {\bf D}_{\rm EnKF} \left({\bf I} - \frac{1}{M} {\bf 1} 
{\bf 1}^{\rm T} \right) + {\bf D}_0
\end{equation}
leads to $\widehat{\bf D}_{\rm EnKF} \in \mathcal{D}_1$.
\end{remark}

Note that the analysis covariance matrix (\ref{cma}) can be equivalently 
written in the form
\begin{equation} \label{cma2}
\widehat{\bf P}^{\rm a} = \frac{1}{M} {\bf Z}^{\rm f} ({\bf D} - {\bf w}{\bf 1}^{\rm T})
({\bf D} - {\bf w}{\bf 1}^{\rm T})^{\rm T} ({\bf Z}^{\rm f})^{\rm T}
\end{equation}
for any ${\bf D}\in\mathcal{D}_1$. In order to achieve second-order accuracy, 
(\ref{cma2}) has to be equal to the importance sampling estimate of the 
posterior covariance matrix (\ref{cme}) which can now be expressed in the following form
\begin{equation} \label{cme2}
{\bf P}^{\rm a} = {\bf Z}^{\rm f} ({\bf W} - {\bf w}{\bf w}^{\rm T} ) ({\bf Z}^{\rm f})^{\rm T}.
\end{equation}
The class of second-order accurate LETFs, considered in this paper, is now characterized by the set 
\begin{equation}
\mathcal{D}_2=\{{\bf D} \in\mathcal{D}_1|\,
({\bf D} - {\bf w}{\bf 1}^{\rm T})
({\bf D} - {\bf w}{\bf 1}^{\rm T})^{\rm T} = {\bf W} - {\bf w}{\bf w}^{\rm T} \,\}.
\end{equation}

\begin{remark}
There is an important subclass $\mathcal{D}^{+}_1\subset\mathcal{D}_1$ that satisfies the additional constraint $d_{ij}\ge0$, i.e.
\begin{equation}
\mathcal{D}_1^{+}=\{{\bf D} \in\mathcal{D}_1|\,d_{ij} \ge 0\,\,\mbox{for all}\,\,i,j=1,\ldots,M\}.
\end{equation}
Then ${\bf D}\in\mathcal{D}^{+}_1$ are left stochastic matrices and thus can be interpreted as resampling schemes that produce realizations ${\bf z}^a_j$ with respect to the transition probabilities in column $j$ in ${\bf D}$ for $j\in{1,\dots,M}$. However, if such a stochastic matrix is used deterministically to produce an analysis ensemble, such as in the ETPF, then the 
particles ${\bf z}^a_j$ are associated with the expected value of the random variable induced by each column $j$ of ${\bf D}$. Consider, for example, the simple transformation matrix 
${\bf D}_0 \in \mathcal{D}^{+}_1$ given in (\ref{simple}). 
In this case, ${\bf z}_j^{\rm a} = \overline{\bf z}^{\rm a}$ for all $j \in \{1,\ldots,M\}$ and the 
implied analysis covariance matrix (\ref{cma2}) becomes identical to zero, which is clearly 
undesirable, and ${\bf D}_0 \notin \mathcal{D}_2$. The ETPF is designed such that this effect 
is minimized and vanishes asymptotically as $M\to \infty$ \cite{sr:reich13,sr:reichcotter15}. 
More broadly speaking, one has $\mathcal{D}_1^{+} \cap \mathcal{D}_2 = \emptyset$ 
generically.
\end{remark}

We now propose a general methodology 
of how to turn a transformation matrix ${\bf D} \in \mathcal{D}_1$ into a transformation matrix $\widehat{\bf D}\in \mathcal{D}_2$.
We start from the \textit{ansatz}
\begin{equation} \label{modupdate}
\widehat{\bf D} = {\bf D} + {\bf \Delta}
\end{equation}
with ${\bf D}\in\mathcal{D}_1$, ${\bf \Delta} \in \mathbb{R}^{M\times M}$ such that ${\bf \Delta}{\bf 1} = {\bf 0}$, ${\bf \Delta}^{\rm T}{\bf 1} = {\bf 0}$, and
${\bf P}^{\rm a} = \widehat{\bf P}^{\rm a}$
with
\begin{equation}
\widehat{\bf P}^{\rm a} = \frac{1}{M} {\bf Z}^{\rm f} (\widehat{\bf D} - {\bf w}{\bf 1}^{\rm T})
(\widehat{\bf D} - {\bf w}{\bf 1}^{\rm T})^{\rm T} ({\bf Z}^{\rm f})^{\rm T}\,.
\end{equation}
The condition
\begin{equation}
{\bf 0} = {\bf P}^{\rm a} - \widehat{\bf P}^{\rm a} = 
{\bf Z}^{\rm f} \left\{ ({\bf W} - {\bf w}{\bf w}^{\rm T} ) 
- \frac{1}{M}(\widehat{\bf D} - {\bf w}{\bf 1}^{\rm T})
(\widehat{\bf D} - {\bf w}{\bf 1}^{\rm T})^{\rm T} \right\}
({\bf Z}^{\rm f})^{\rm T} ,
\end{equation}
together with (\ref{modupdate}) lead to the following quadratic equation in the correction ${\bf \Delta}$:
\begin{equation} \label{riccati}
M ({\bf W} - {\bf w}{\bf w}^{\rm T} ) 
-({\bf D} - {\bf w}{\bf 1}^{\rm T})
({\bf D} - {\bf w}{\bf 1}^{\rm T})^{\rm T}  =
({\bf D}-{\bf w}{\bf 1}^{\rm T}) {\bf \Delta}^{\rm T} + {\bf \Delta} ({\bf D} - {\bf w}{\bf 1}^{\rm T})^{\rm T}
 + {\bf \Delta}{\bf \Delta}^{\rm T}.
\end{equation}
If we also choose ${\bf \Delta}$ to be symmetric, then the special case (\ref{simple}) leads to 
\begin{equation}
M({\bf W} - {\bf w}{\bf w}^{\rm T})={\bf \Delta}{\bf \Delta}
\end{equation}
and a solution of (\ref{riccati}) is simply given by the symmetric square root
\begin{equation}
{\bf \Delta} = \sqrt{M} ({\bf W} - {\bf w}{\bf w}^{\rm T})^{1/2},
\end{equation}
which recovers the NETF \cite{sr:toedter15,sr:xiong06}. 
Note that ${\bf \Delta} {\bf Q}$ with ${\bf Q}$ an $M\times M$ orthogonal matrix 
such that ${\bf Q}{\bf 1} = {\bf 1}$ also provide a solution to (\ref{riccati}) 
if ${\bf D}={\bf w}{\bf 1}^{\rm T}$.\footnote{The NETF, as proposed in \cite{sr:toedter15}, uses randomly chosen orthogonal matrices which satisfy 
${\bf Q}{\bf 1}  = {\bf 1}$ while the NETF of 
\cite{sr:xiong06} is based on a non-symmetric square root of ${\bf W} - {\bf w}{\bf w}^{\rm T}$.} 
The following lemma states how to choose the orthogonal matrix ${\bf Q}$ in an optimal way.

\begin{lemma}\label{lemma1}
Let ${\bf \Delta}$ be any $M\times M$ matrix such that (i) ${\bf \Delta} {\bf 1} = {\bf 0}$ and (ii)
\begin{equation}
\frac{1}{M} {\bf \Delta} {\bf \Delta}^{\rm T} = {\bf W} - {\bf w}{\bf w}^{\rm T}
\end{equation}
and let us assume that $M\le N_z + 1$.
Define the $M\times M$ orthogonal matrix
\begin{equation} \label{optimal}
 {\bf Q}_{\rm opt}:= {\bf U}_{\rm opt} {\bf V}_{\rm opt}^{\rm T}
\end{equation}
with the two $M\times M$ orthogonal matrices ${\bf U}_{\rm opt}$ and ${\bf V}_{\rm opt}$ given by the singular value decomposition of the 
$M\times M$ matrix  
\begin{equation}\label{perturbations}
{\bf S} = {\bf \Delta} (\widehat{\bf Z}^{\rm f})^{\rm T} \widehat{\bf Z}^{\rm f}, \qquad
\widehat{\bf Z}^{\rm f} := {\bf Z}^{\rm f} - \frac{1}{M} {\bf Z}^{\rm f}{\bf 1} {\bf 1}^{\rm T},
\end{equation}
i.e.~${\bf S} = {\bf U}_{\rm opt} {\bf \Lambda}_{\rm opt} {\bf V}_{\rm opt}^{\rm T}$. 
Then the transformation matrix
\begin{equation} \label{optimalD1}
{\bf D}_{\rm opt} = {\bf w} {\bf 1}^{\rm T} + {\bf \Delta} {\bf Q}_{\rm opt}
\end{equation}
results in a second-order accurate LETF, which minimizes
\begin{equation} \label{cost2}
\widehat{J} ({\bf D}) = \frac{1}{M}  \sum_{i=1}^M \|{\bf z}_i^{\rm a} - {\bf z}_i^{\rm f} \|^2
\end{equation}
over all second-order accurate transformation matrices ${\bf D}\in\mathcal{D}_2$. 
\end{lemma}

\noindent
{\it Proof.}
Since $\widehat{\bf Z}^{\rm f} {\bf 1} = {\bf 0}$, the matrix ${\bf S}$ also satisfies 
${\bf S}{\bf 1} = {\bf 0}$ in addition to ${\bf S}^{\rm T}{\bf 1} = {\bf 0}$, which
implies that ${\bf Q}_{\rm opt}{\bf 1} = {\bf 1}$ and (\ref{optimalD1}) is second-order accurate. Also note that
\begin{eqnarray}
  \left( {\bf \Delta} (\widehat{\bf Z}^{\rm f})^{\rm T}
\widehat{\bf Z}^{\rm f} (\widehat{\bf Z}^{\rm f})^{\rm T} \widehat{\bf Z}^{\rm f} {\bf \Delta}\right)^{-1/2} 
{\bf \Delta} (\widehat{\bf Z}^{\rm f})^{\rm T} \widehat{\bf Z}^{\rm f} &=&
\left( {\bf S} {\bf S}^{\rm T} \right)^{-1/2} {\bf S}\\
&=& \left({\bf U}_{\rm opt} {\bf \Lambda}_{\rm opt}^{-1} {\bf U}^{\rm T}_{\rm opt} \right)
{\bf U}_{\rm opt} {\bf \Lambda}_{\rm opt}{\bf V}^{\rm T}_{\rm opt} 
\end{eqnarray} 
which has been shown in \cite{sr:olkin82} to minimize (\ref{cost2}) for given forecast and analysis means and 
covariance matrices and the optimality of ${\bf Q}_{\rm opt} = {\bf U}_{\rm opt} {\bf V}_{\rm opt}^{\rm T}$ follows. See also \cite{sr:reichcotter15}. \hfill $\Box$ \\

A couple of comments should be made on the requirement of $M\le N_z + 1$ in Lemma 
\ref{lemma1}. First, if the number of samples, $M$, exceeds the dimensions of 
state space, $N_z$, then it is 
computationally preferable to implement the optimal transformation in the form
\begin{equation}
{\bf z}_i^{\rm a} = \overline{\bf z}^{\rm a} + {\bf T} ({\bf z}_i^{\rm f} - \overline{\bf z}^{\rm f}),
\end{equation}
where ${\bf T} \in \mathbb{R}^{N_z\times N_z}$ is an appropriately defined symmetric matrix \cite{sr:olkin82,sr:reichcotter15}. Second, one could still proceed with (\ref{optimal}) but
should multiply ${\bf Q}_{\rm opt}$ by the projection matrix ${\bf I}- {\bf 1}{\bf 1}^{\rm T}/M$ 
from the right in order to keep the resulting transformation matrix (\ref{optimalD1}) 
mean preserving, i.e., $\mathbf{D}_{\rm opt} {\bf 1}={\bf w}$. 
This additional operation arises from the fact that the matrix ${\bf S}$ 
will have multiple zero singular values.


\section{Continuous-time algebraic Riccati equation} \label{sec:Riccati}

We now return to the general case of a first-order accurate transformation matrix ${\bf D}$. Then
(\ref{riccati}) leads to a continuous-time 
algebraic Riccati equation in the symmetric correction ${\bf \Delta}$. 
More specifically, upon introducing
\begin{equation} \label{AB}
{\bf B} = {\bf D} - {\bf w}{\bf 1}^{\rm T}, \qquad \qquad {\bf A} = M({\bf W}-{\bf w} {\bf w}^{\rm T}) - {\bf B}{\bf B}^{\rm T}
\end{equation}
and assuming that ${\bf \Delta}$ is symmetric, equation
(\ref{riccati}) can be expressed as the continuous-time algebraic Riccati equation
\begin{equation} \label{riccati2}
{\bf A} = {\bf B} {\bf \Delta} + {\bf \Delta}{\bf B}^{\rm T} + {\bf \Delta}
{\bf \Delta}.
\end{equation}
Note that (\ref{riccati2}) arises as the stationary solution of the dynamic Riccati equation
\begin{equation} \label{driccati}
\frac{\rm d}{{\rm d}\tau} {\bf \Delta} = - {\bf B} {\bf \Delta} - {\bf \Delta}{\bf B}^{\rm T} + 
{\bf A} - {\bf \Delta}{\bf \Delta}.
\end{equation}
Since (\ref{driccati}) is controllable \cite{sr:Wonham68}, solutions, ${\bf \Delta}(\tau)$, 
of (\ref{driccati}) with initial condition ${\bf \Delta}(0) = {\bf 0}$
will converge to a solution of (\ref{riccati2}) as $\tau \to \infty$ \cite{sr:Bucy}. Hence numerical time-stepping of (\ref{driccati}) with the explicit Euler method 
for sufficiently many iterations will result in an approximate solution to (\ref{riccati2}). This approach has been used for the numerical results
displayed later in this paper. 

\begin{remark}
Alternatively, (\ref{riccati2}) can be solved by applying the Schur vector approach of \cite{sr:laub79}. 
The Schur vector approach is based on the extended Hamiltonian matrix
\begin{equation} \label{schur}
{\bf H} = \left( \begin{array}{cc} {\bf B}^{\rm T} & {\bf I} \\ {\bf A} & -{\bf B} \end{array} \right)
\end{equation}
and its upper triangular Schur decomposition
\begin{equation} \label{Schurdecomp}
{\bf U}^{\rm T} {\bf H} {\bf U} = \left( \begin{array}{cc} {\bf S}_{11} & {\bf S}_{12} \\ {\bf 0} &
{\bf S}_{22} \end{array} \right)
\end{equation}
with the real part of the spectrum of ${\bf S}_{11}$ being negative and the real parts of the spectrum of
${\bf S}_{22}$ being positive. With the orthogonal matrix ${\bf U}$ partitioned accordingly, 
the solution of (\ref{riccati}) is given by
\begin{equation} \label{correction1}
{\bf \Delta} = {\bf U}_{21} {\bf U}_{11}^{-1}.
\end{equation}
This computational approach requires that the  matrix pair 
$({\bf A}^{1/2},{\bf B})$ is detectable \cite{sr:Wonham68,sr:laub79}. Since this condition
may not always be satisfied for (\ref{riccati2}), we recommend to use (\ref{driccati}) in order to find approximative solutions to (\ref{riccati2}). Alternatively, one could exploit more general Lagrangian
invariant subspace techniques as discussed in \cite{sr:FMX02}.
\end{remark}

We can now either use the ${\bf D}$ from the ETPF and derive a second-order accurate version of the ETPF or we  compute an approximate solution ${\bf D} \in \mathcal{D}_1^{+}$ 
to the optimal transport problem (\ref{cost}) and are still  able to turn it into a second-order 
accurate PF. This aspect will be discussed in more detail in Section \ref{sec:sinkhorn}.


%
\section{Sinkhorn approximation to the optimal transport problem} \label{sec:sinkhorn}
%

The Sinkhorn approximation to the optimal transport problem defined by 
the cost functional (\ref{cost}) and ${\bf D}\in \mathcal{D}_1^{+}$ is provided by
the regularised cost functional
\begin{equation} \label{costSH}
\mathbf{D}_{\rm min}(\lambda )=\arg \min J_{\rm SH}({\bf D}) = \sum_{i,j=1}^M \left\{ 
d_{ij} \|{\bf z}_i^{\rm f}-{\bf z}_j^{\rm f}\|^2 + \frac{1}{\lambda} d_{ij} \ln \frac{d_{ij}}{d^0_{ij}} \right\}
\end{equation}
where $\lambda > 0$ is a regularization parameter and $d^0_{ij}$ are the entries of $\mathbf{D}_0$ defined in (\ref{simple}). Each parameter $\lambda$ is associated with a specific $\mathbf{D}_{\rm min}(\lambda ) \in \mathcal{D}_1^{+}$ and $\lambda \to \infty$ leads back to the original cost function (\ref{cost}). While, one the other hand, the choice $\lambda \to 0$ leads to (\ref{simple}) as the unique minimizer. This follows from the fact that  the regularization term 
in (\ref{costSH}) is minimal for $d_{ij}=d^0_{ij} $, i.e., $\lim_{\lambda\rightarrow 0}\mathbf{D}_{\rm min}(\lambda ) = \mathbf{D}_0$.

\begin{remark}
After determining ${\bf D}_{\rm min}(\lambda)$ it is possible to add an appropriate corresponding second-order correction term ${\bf \Delta}(\lambda)$ which, depending on $\lambda$, leads to different second-order accurate particle filters, e.g., $\lambda \rightarrow 0$ leads to the NETF and $\lambda \rightarrow
\infty$ to the second-order corrected ETPF. In other words, varying $\lambda$ allows one to naturally 
bridge between the NETF and the second-order corrected ETPF.
\end{remark}

There exists a straightforward iterative method for finding the minimizer
of (\ref{costSH}). First one notes that the minimizer is of the form
\begin{equation}
{\bf D}_{\rm min}(\lambda) = \mbox{diag}({\bf u}) \,{\bf K} \,\mbox{diag}({\bf v}),
\end{equation}
where ${\bf u} \in \mathbb{R}^{M\times 1}$ and ${\bf v}\in \mathbb{R}^{M\times 1}$ 
are two non-negative  vectors  and ${\bf K}$ has entries
\begin{equation} \label{entriesK}
k_{ij} = e^{-\lambda \|{\bf z}_i^{\rm f} - {\bf z}_j^{\rm f}\|^2 }.
\end{equation}
The unknown vectors ${\bf u}$ and ${\bf v}$ can be computed by Sinkhorn's fixed point iteration
\begin{equation} \label{sinkhorn}
\{ M w_i/({\bf K} {\bf v})_i \} \to {\bf u} ,\qquad \{1/({\bf K} {\bf u})_i\} \to {\bf v}.
\end{equation}
The Sinkhorn approximation requires ${\cal O}(M^2)$ operations. See \cite{sr:cuturi13} 
for an efficient implementation and additional details.

Let us denote the iterates of ${\bf u}$ and ${\bf v}$ by ${\bf u}^l$ and ${\bf v}^l$, respectively, 
where we always update ${\bf u}$ first according to the formula to the left in (\ref{sinkhorn}). Then the associated
\begin{equation} \label{sinkhornD}
{\bf D}^l =  \mbox{diag}({\bf u}^l) \,{\bf K} \,\mbox{diag}({\bf v}^l)
\end{equation}
satisfies $({\bf D}^l)^{\rm T}{\bf 1} = {\bf 1}$ and the weights
\begin{equation} \label{sinkhornw}
{\bf w}^l = \frac{1}{M} {\bf D}^l {\bf 1}
\end{equation}
converge to ${\bf w}$ as $l \to \infty$. If we stop the iteration at an index $l_\ast$, then we define
the associated transformation matrix by
\begin{equation} \label{sinkhornfinal}
{\bf D} = {\bf D}^{l_\ast} - ({\bf w}^{l_\ast}  + {\bf w}) {\bf 1}^{\rm T}.
\end{equation}
The index $l_\ast$ can be determined by the condition
\begin{equation} \label{sinkhornstop}
\| {\bf w}^{l_\ast}- {\bf w}\| \le \varepsilon 
\end{equation}
for sufficiently small $\varepsilon >0$, e.g.~$\varepsilon = 10^{-8}$.

%
\section{Algorithmic summary} \label{sec:AS}
%

We summarise the key steps of the second-order accurate ETPF implementation based upon
the Sinkhorn approximation to the optimal transport problem. The Sinkhorn approximation
can, of course, be replaced by any available direct solver for the optimal transport problem.

We assume that a set of forecast ensemble members, ${\bf Z}^{\rm f}$, and
a vector of importance weights, ${\bf w}$, are given. Then the following steps are performed:

\begin{itemize}
\item[(i)] Select a regularization parameter $\lambda > 0$ for the 
Sinkhorn approximation to the optimal transport algorithm. Compute the matrix
${\bf K}$ according to (\ref{entriesK}). Normalize the entries of ${\bf K}$ such that all
entries satisfy $-\lambda^{-1} \ln k_{ij} \le 1$. Recursively compute vectors ${\bf u}^l$ and ${\bf v}^l$
according to the update formula (\ref{sinkhorn}). Start with ${\bf v}^0 = {\bf 1}$. 
Iterate till the transformation matrix (\ref{sinkhornD}) and its associated weight vector
(\ref{sinkhornw}) satisfy (\ref{sinkhornstop}). Note that (\ref{sinkhornD}) should satisfy
${\bf 1}^{\rm T}{\bf D}^l = {\bf 1}^{\rm T}$ in each iteration. We used 
$\varepsilon = 10^{-8}$ in our experiments. One finally obtains the transform matrix
${\bf D}$ using (\ref{sinkhornfinal}).
\item[(ii)] Solve the Riccati equation (\ref{riccati2}) for the correction ${\bf \Delta}$ by solving
the dynamic Riccati equation (\ref{driccati}) with the explicit Euler method, 
step-size $\Delta \tau = 0.1$, and initial condition ${\bf \Delta}(0) = {\bf 0}$. The iteration is stopped
whenever
\begin{equation}
\|{\bf \Delta}((k+1)\Delta \tau)-{\bf \Delta}(k\Delta \tau)\|_\infty \le 10^{-3}
\end{equation}
and we set ${\bf \Delta} = {\bf \Delta}((k+1)\Delta \tau)$.
\item[(iii)] The analysis ensemble is given by
\begin{equation}
{\bf Z}^{\rm a} = {\bf Z}^{\rm f} \widehat{\bf D} = {\bf Z}^{\rm f}({\bf D} + {\bf \Delta}).
\end{equation}
\end{itemize}

We mention that the proposed second-order accurate ETPF can be
used instead of the standard ETPF in a hybrid filter, as described in
\cite{sr:CRR15}, and, when applied to spatially extended system, can
also be used with localization. More specifically, a hybrid filter is based on
factorizing the likelihood (\ref{likelihood}) into
\begin{equation}
\pi({\bf y}|{\bf z}) = \pi({\bf y}|{\bf z})^\alpha \times \pi({\bf y}|{\bf z})^{1-\alpha}
\end{equation}
and applying different filters to each of the two factors. $R$-localization, on the other hand, 
leads to different transformation matrices ${\bf D}(x_k)$ at each grid point $x_k$ of
the computational domain. See \cite{sr:reich15,sr:reichcotter15} for further details.

%
\section{Numerical examples} \label{sec:NE}
%


We now demonstrate the numerical behavior of the proposed second-order accurate ETPF
as summarized in Section \ref{sec:AS}.  The first two experiments are based on the Lorenz-63 
and Lorenz-96 models, respectively, and its data assimilation setting of
\cite{sr:CRR15}. We finally apply the second-order accurate filters to parameter estimation of
the scene-viewing model {\it SceneWalk} \cite{sr:Engbert2015}.

\subsection{Lorenz-63} \label{sec:num1}

\begin{figure}
\begin{center}
\includegraphics[width=0.47\textwidth]{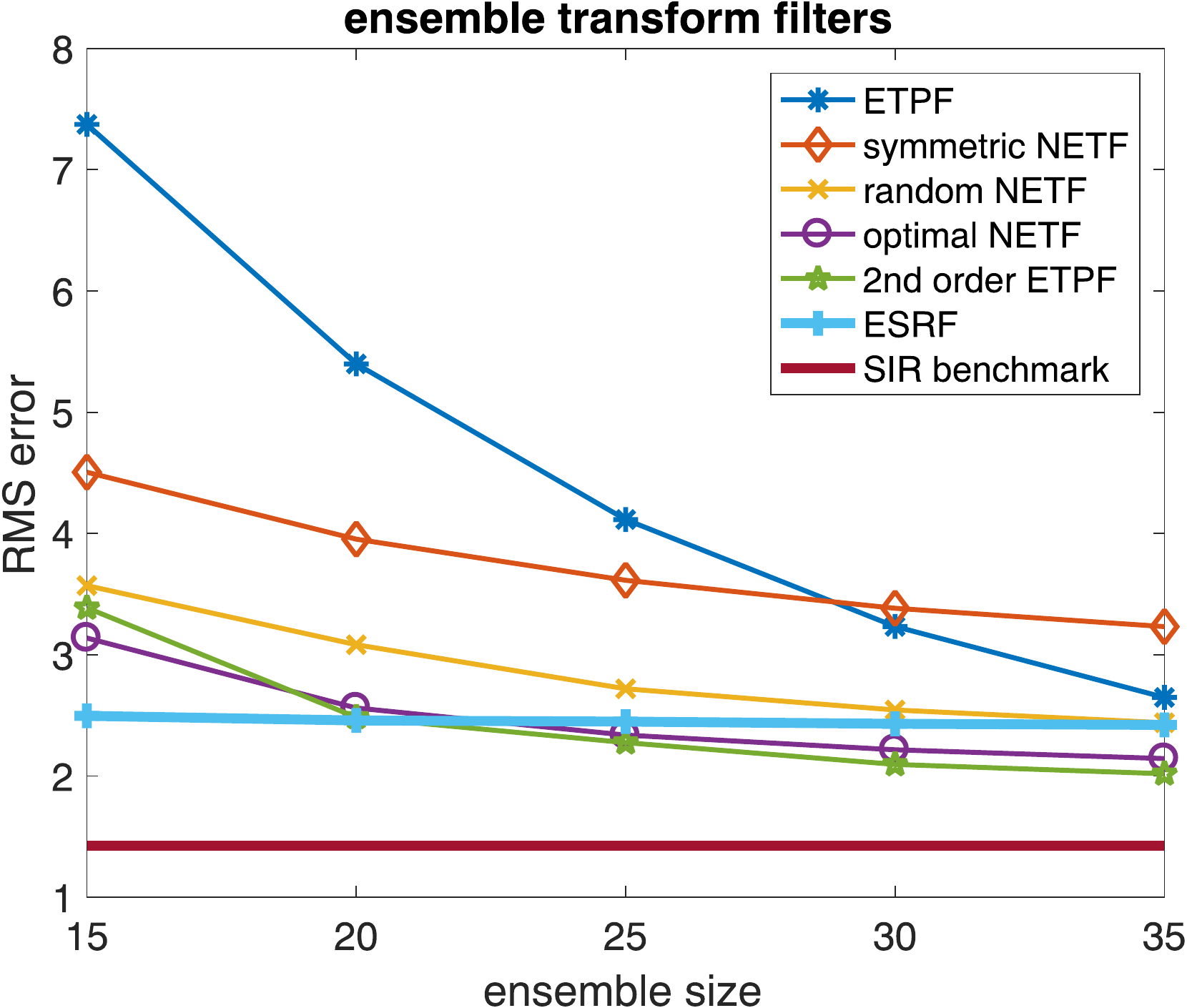} $\qquad$
\includegraphics[width=0.47\textwidth]{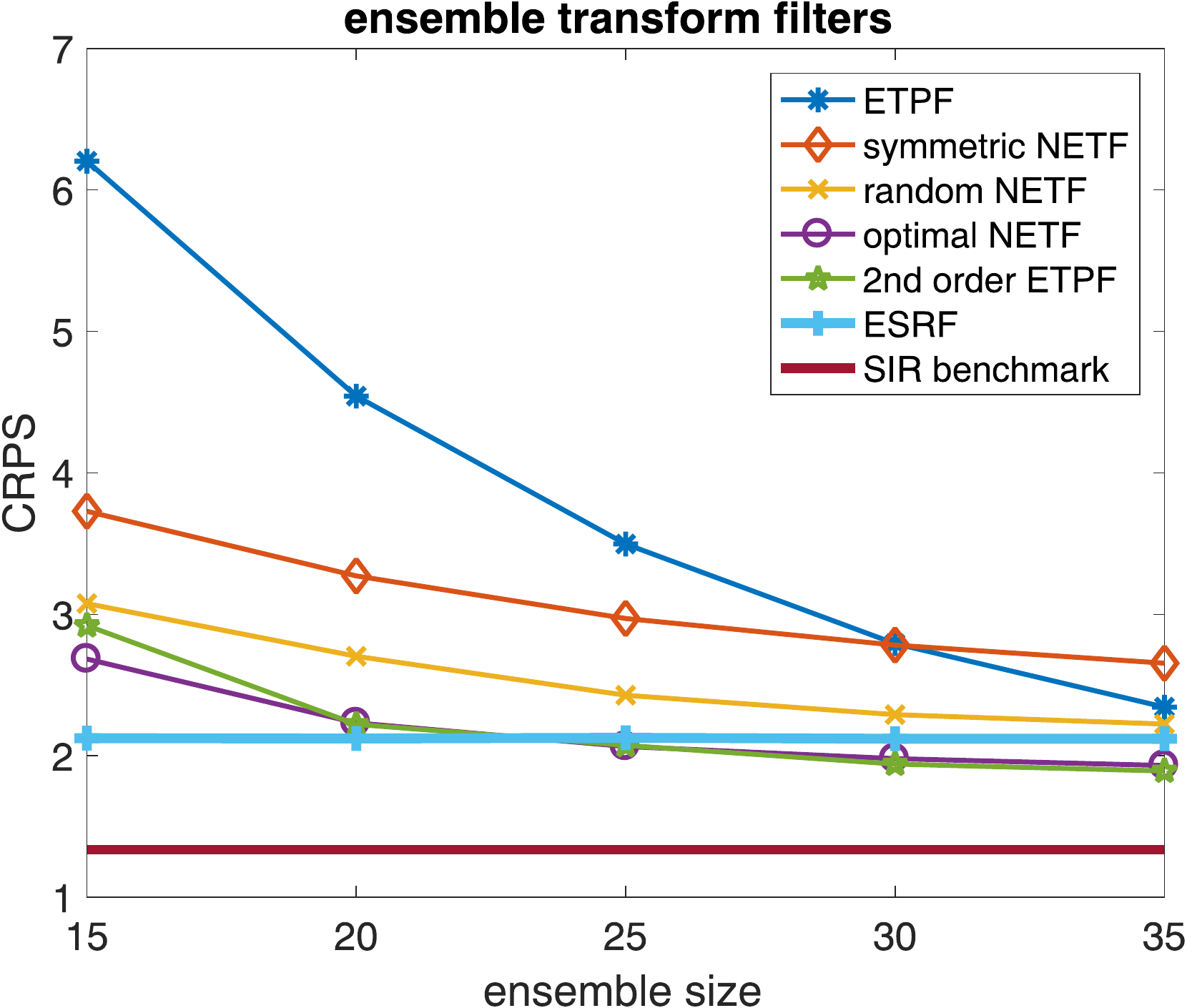} \end{center}
\caption{RMS errors (left panel) and CRPS (right panel) 
for various second-order accurate LETFs compared to the
ETPF and the ESRF as a function of the ensemble size, $M$, for the Lorenz-63 model. We also provide the RMS error
and the CRPS obtained from a standard particle filter with resampling and $M = 1000$ ensemble members.}
\label{figure2.0}
\end{figure}

We use the chaotic Lorenz-63 system \cite{sr:lorenz63} with the standard parameter setting 
$\sigma = 10$, $\rho = 28$, and $\beta = 8/3$, and observe the first component
of the three dimensional system in observation intervals of $\Delta t_{\rm obs}=0.12$ with observation 
error variance $R=8$. A total of $K=500,000$ assimilation steps are performed. 
Since the model dynamics is deterministic, particle rejuvenation
\begin{equation} \label{rejuvenation}
{\bf z}_j^{\rm a} \to {\bf z}_j^{\rm a} + \sum_{i=1}^M ({\bf z}_i^{\rm f} - \bar {\bf z}^{\rm f}) \frac{\beta \xi_{ij}}{\sqrt{M-1}}
\end{equation}
is applied with $\beta = 0.2$ and $\xi_{ij}$ independent and identically distributed Gaussian random variables with
mean zero and variance one. Simulations with $\beta = 0.15$ and $\beta = 0.25$ gave similar results to those reported here. This data assimilation setting
has already been used in \cite{sr:CRR15} and \cite{sr:reich15} 
since it leads to non-Gaussian forecast and analysis distributions and a particle filter is able to outperform EnKFs in the limit of large ensemble sizes. 

A comparison between a standard particle
filter with resampling, the EnKF, and the ETPF can be found in \cite{sr:reich15}.
Here we are, however, interested in the performance of second-order accurate filters for small ensemble sizes in the range $M\in \{15,20,\ldots,35\}$. See Figure \ref{figure2.0} for the 
resulting time-averaged RMS errors. It can be clearly seen that the 
second-order corrected ETPF and the NETF with optimally chosen rotation matrix  leads to 
the smallest RMS errors for $M\ge 25$, while the 
standard ensemble square root filter (ESRF) \cite{sr:evensen} is optimal for smaller ensemble sizes. 
It can be seen that the standard ETPF is not competitive except for $M=35$. 
The same findings apply for the continuous ranked probability score (CRPS)
\cite{sr:broecker12}, which we computed for the observed component of the Lorenz-63 system.
The results can be found in Figure \ref{figure2.0}. 

We also display the RMS errors for implementations of the NETF with  randomly chosen orthogonal matrices, ${\bf Q}$, as suggested by \cite{sr:toedter15}, and with ${\bf Q} = {\bf I}$ 
in Figure \ref{figure2.0}. It can be seen that both choices lead to substantially increased RMS
errors.

\begin{figure}
\begin{center}
\includegraphics[width=0.47\textwidth]{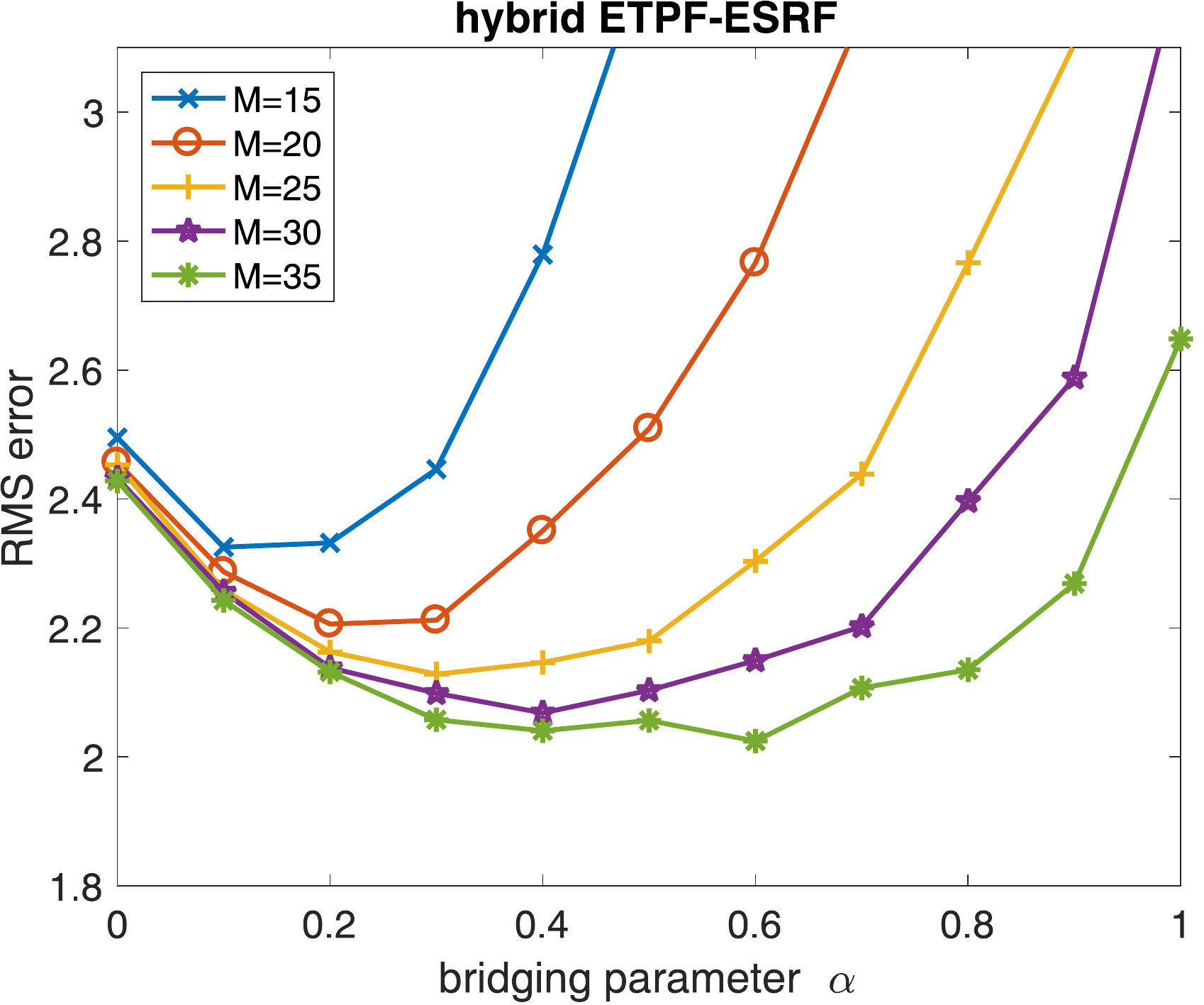} $\qquad$
\includegraphics[width=0.47\textwidth]{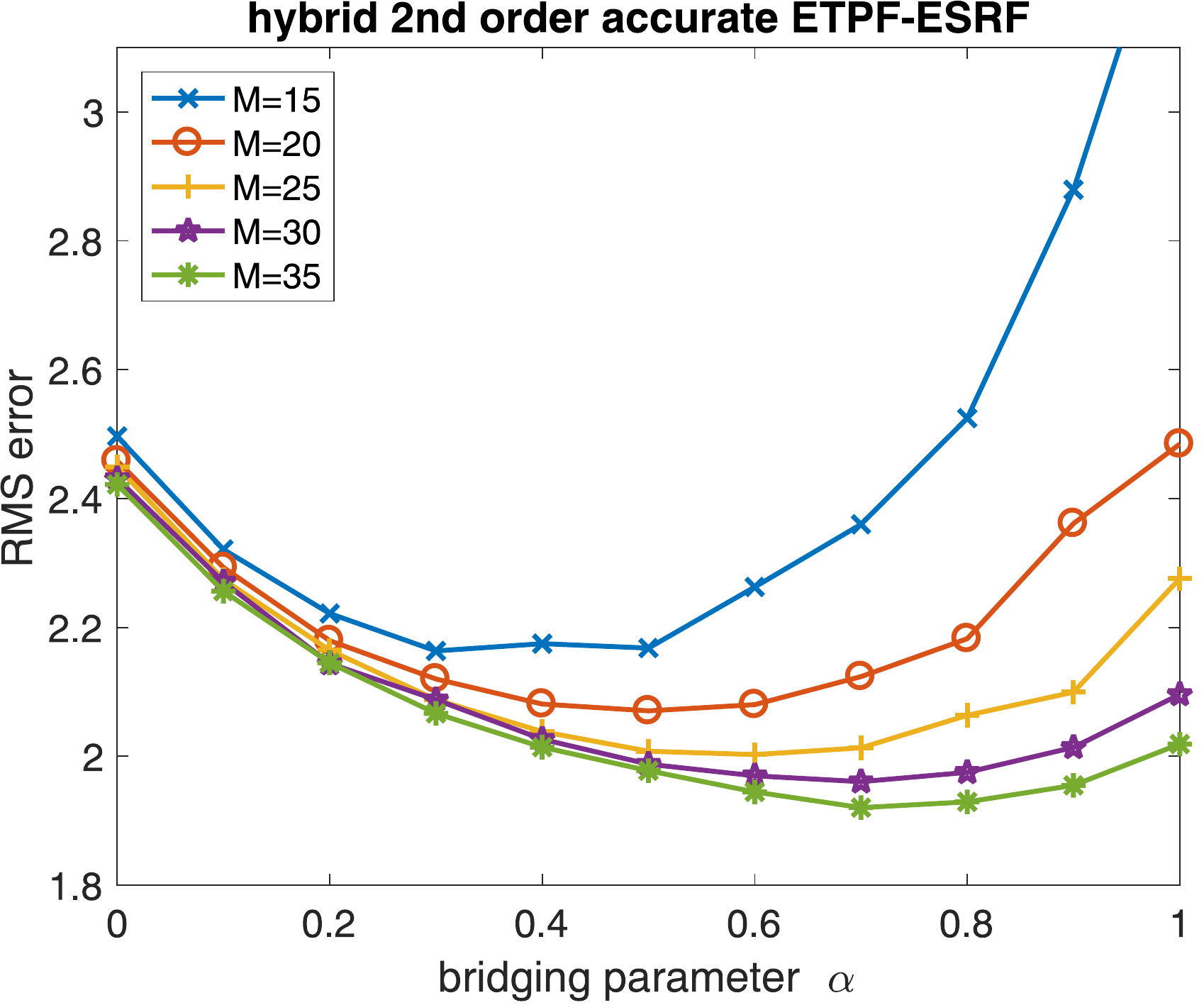}
\end{center}
\caption{Hybrid filter with standard EPTF (left panel) and second-order accurate ETPF (right
panel) applied to the Lorenz-63 model. Time-averaged RMS errors are displayed as a function of the bridging parameter
$\alpha$. Please note that $\alpha=0$ corresponds to the standard ESRF, while $\alpha=1$ corresponds to the  ETPF and the second-order corrected ETPF, respectively. }
\label{figure2.1}
\end{figure}

We now test the second-order accurate transform filters
within the hybrid filter framework proposed in \cite{sr:CRR15}. More specifically,
the hybrid filter of \cite{sr:CRR15} with a second-order accurate transform filter  applied first is implemented for ensemble sizes varying between $M=15$ and $M=35$. The bridging parameter, $\alpha$, of the hybrid filter approach is chosen
such that $\alpha = 0$ corresponds to the standard ESRF 
while $\alpha = 1$ leads to a purely second-order accurate ETPF.
We perform experiments for fixed bridging parameters $\alpha \in \{0,0.1,0.2,\ldots,0.9,1.0\}$
and compare the resulting RMS errors to those from a hybrid method based on the standard
ETPF in Figure \ref{figure2.1}. The improvement achieved by the second-order correction is
clearly visible. In both cases, the ETPF has been implemented using a direct solver
for the underlying optimal transport problem.  

\begin{figure}
\begin{center}
\includegraphics[width=0.47\textwidth]{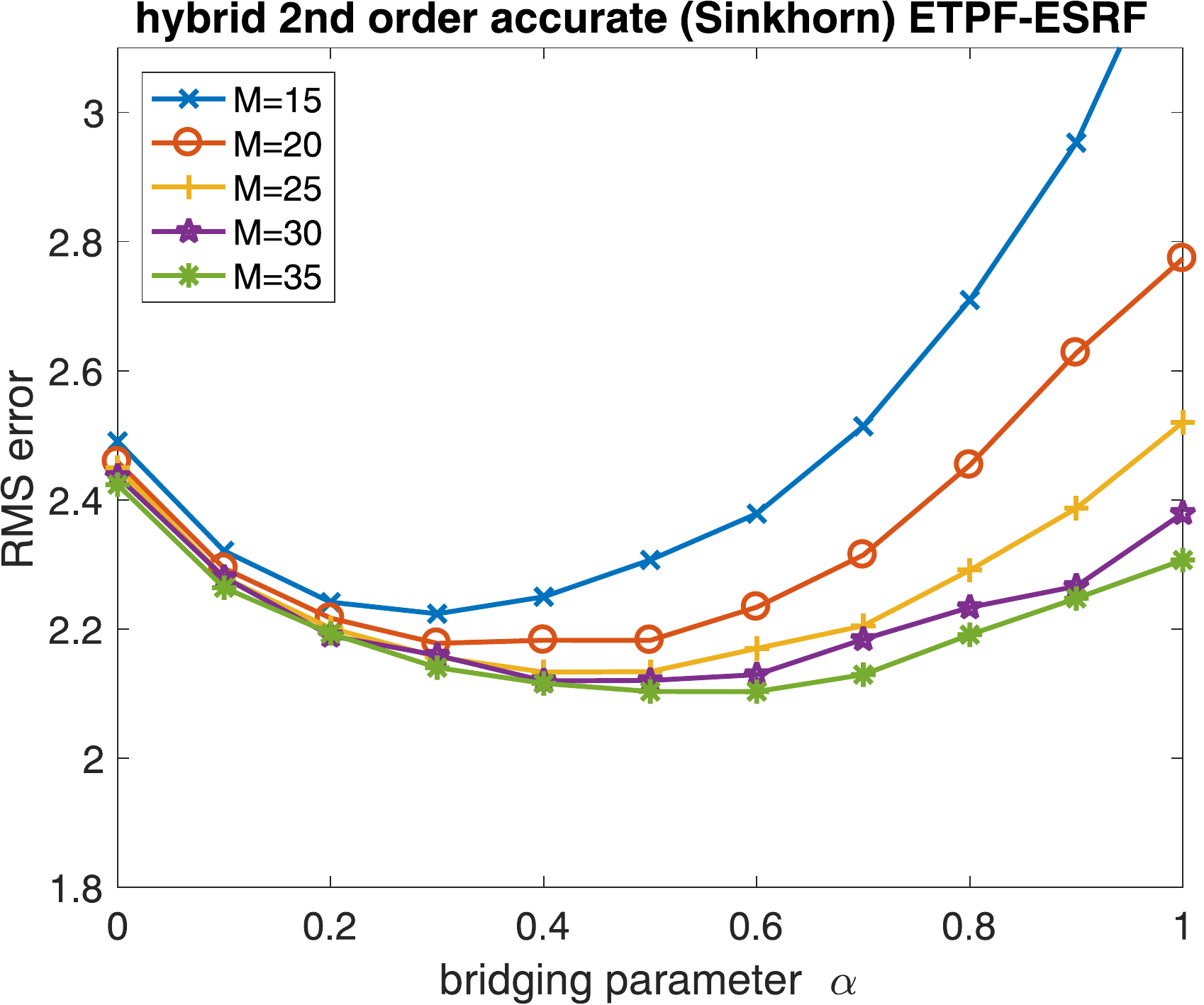} $\qquad$
\includegraphics[width=0.47\textwidth]{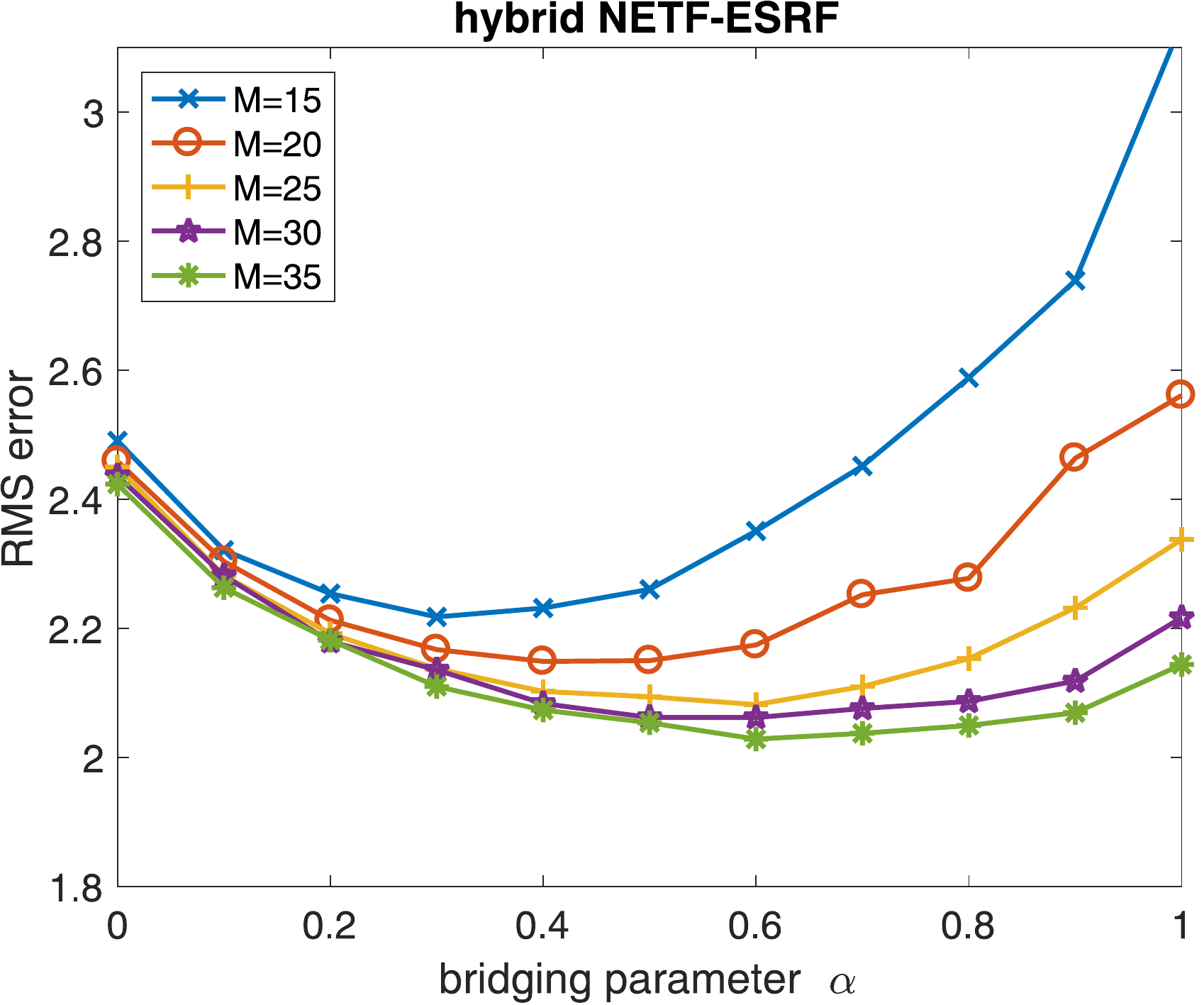}
\end{center}
\caption{Second-order hybrid ETPF-ESRF  with the optimal transport problem
solved by the Sinkhorn approximation with $\lambda = 10$ (left panel) and 
hybrid NETF-ESRF with the  the orthogonal matrix ${\bf Q}$ as defined in (\ref{optimal}) (right 
panel) applied to the Lorenz-63 model. Time-averaged RMS errors are displayed as a function of the bridging parameter
$\alpha$.}
\label{figure2.2}
\end{figure}

We next replace the direct solver for the optimal transport problem by the Sinkhorn approximation
with regularization parameters $\lambda = 10$ and $\lambda = 40$. 
The RMS errors for the resulting hybrid filter with $\lambda = 10$ 
can be found in Figure \ref{figure2.2}, while $\lambda = 40$ leads to RMS errors which
are very close to those displayed in the right panel of Figure \ref{figure2.1}, which are based on
a direct solver for the optimal transport problem. 

We also  implement the hybrid filter of \cite{sr:CRR15} 
with the ETPF being replaced by the 
second-order accurate NETF with the rotation matrix, ${\bf Q}$, defined as in (\ref{optimal}).
We denote this hybrid filter by NETF-ESRF.  The numerical results can also be found in
Figure \ref{figure2.2}. Overall, we find that a second-order corrected hybrid ETPF-ESRF 
and the NETF-ESRF with optimally chosen rotation matrix perform quite comparable in
terms of their RMS errors. The same holds true for the associated CRPS (not displayed).


\subsection{Lorenz-96} \label{sec:num2}

We now implement the spatially-extended Lorenz-96 system \cite{sr:lorenz96} with the standard parameter setting of $p=40$ grid points and forcing $F=8$.
We observe every second grid point in observation intervals of $\Delta t_{\rm obs}=0.11$ with observation error variance $R=8$. 
A total of $K=50,000$ assimilation steps are performed. 
We apply localization \cite{sr:reichcotter15} 
with the localization radius $r_{\rm loc}$ set equal to four grid points and compute
separate transformation matrices $\mathbf{D}(x_k)$ for each grid point $x_k = k$, 
$k=1,\ldots,40$. Localization is necessary for
this test problem as the ensemble sizes, $M \in \{20,25,30\}$, are smaller 
than the number of grid points, $p=40$. This specific DA setting has already been 
used in \cite{sr:reich15} and \cite{sr:CRR15}.

\begin{figure}
\begin{center}
\includegraphics[width=0.47\textwidth]{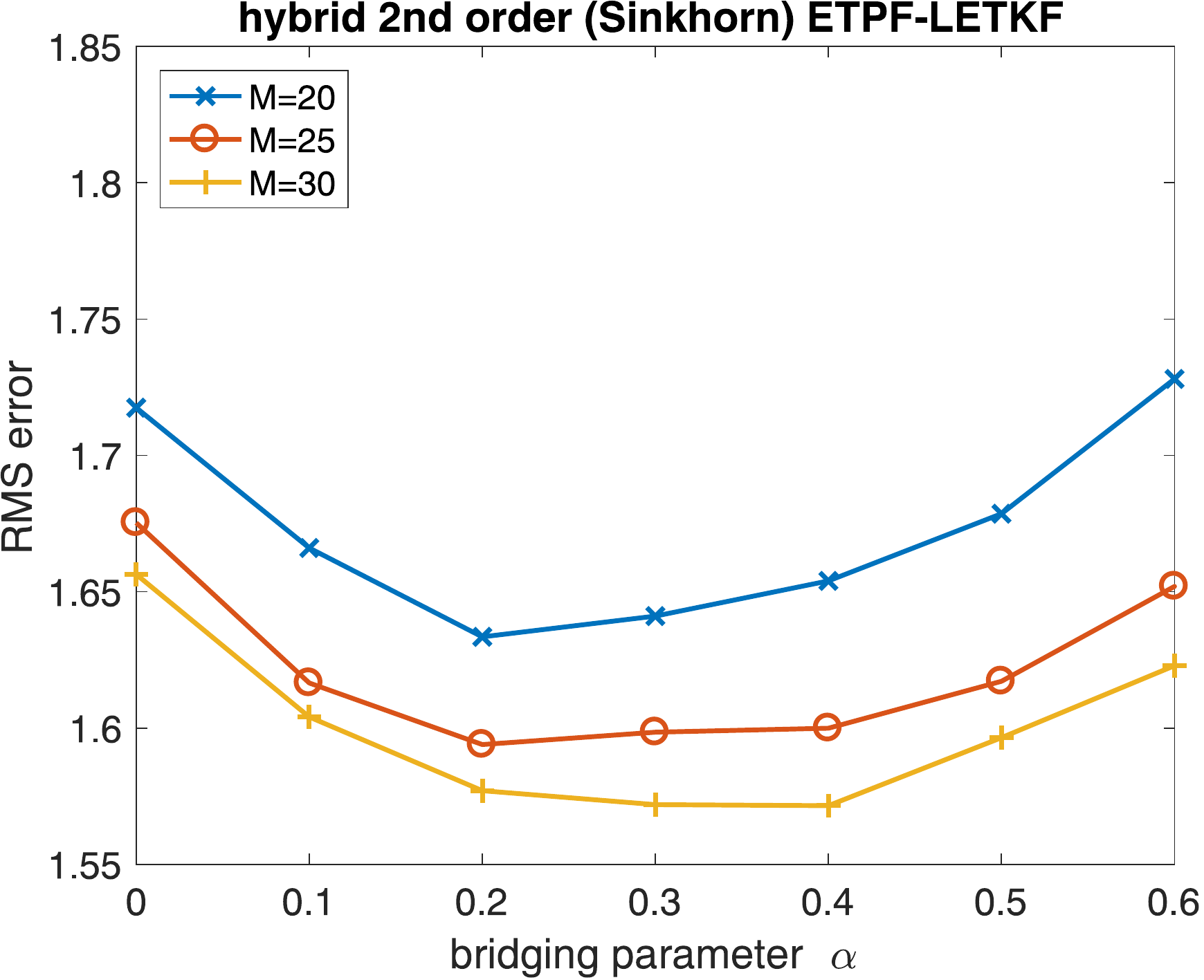} $\qquad$
\includegraphics[width=0.47\textwidth]{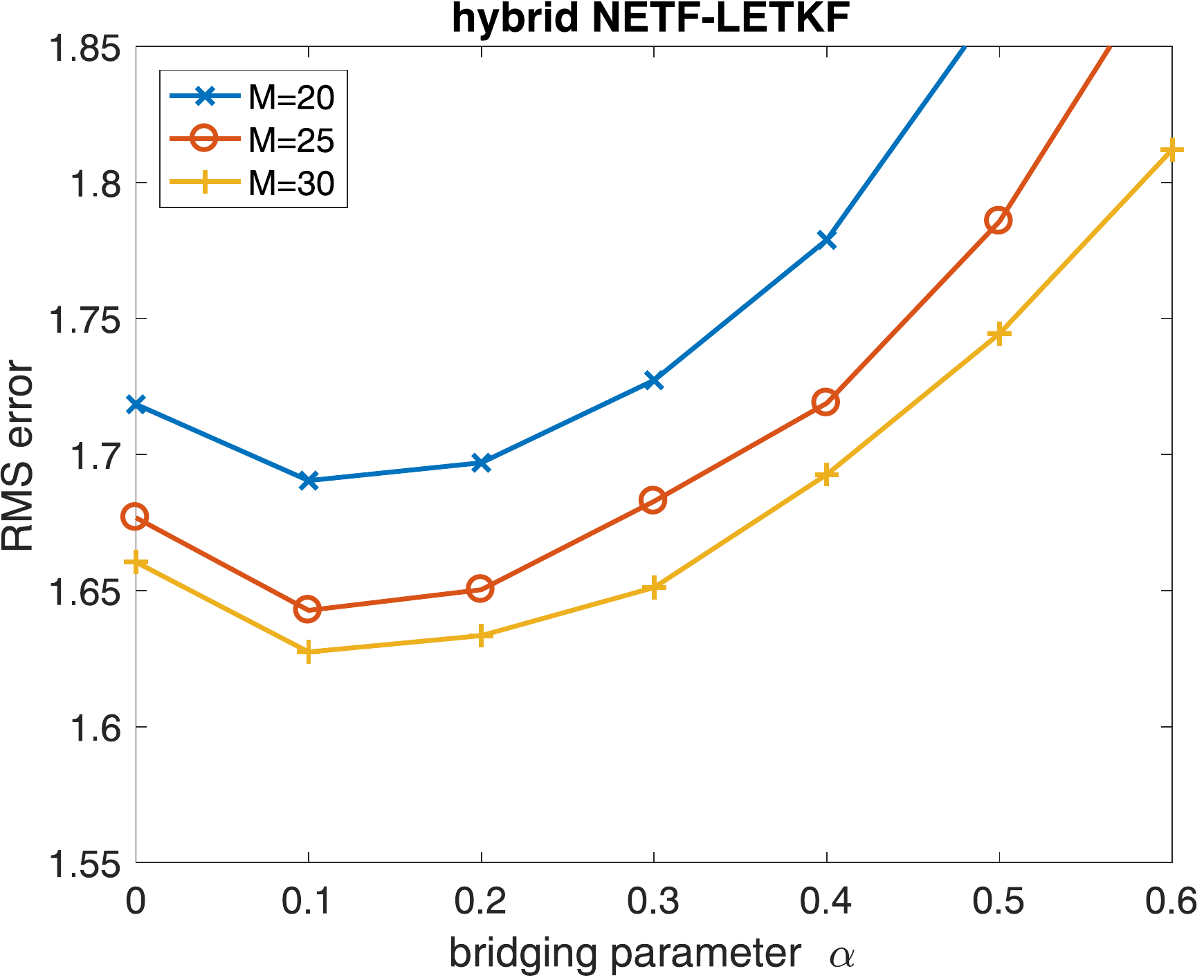}
\end{center}
\caption{Second-order hybrid ETPF-LETKF  with the optimal transport problem
solved by the Sinkhorn approximation with $\lambda = 10$ (left panel) and 
hybrid NETF-LETKF with the  the orthogonal matrix ${\bf Q}$ at each grid point defined as 
in (\ref{optimal}) (right  panel) applied to the Lorenz-96 model. Time-averaged RMS errors are displayed as a function of 
the bridging parameter $\alpha$. The choice $\alpha = 0$ corresponds to the LETKF}
\label{figure4.1}
\end{figure}

We compare two hybrid methods based on a combination of second-order accurate filters
and the local ensemble transform Kalman filter (LETKF) \cite{sr:hunt07}. 
All filters use $R$-localization \cite{sr:hunt07,sr:reichcotter15} and the transportation cost 
at each grid point $x_k = k$, $k=1,\ldots,40$, is given by
\begin{equation}
J({\bf D}(x_k)) = \sum_{i,j=1}^M d_{ij}(x_k) | u_i^{\rm f}(x_k) - u_j^{\rm f}(x_k)|^2,
\end{equation}
where $u_i^{\rm f}(x_k) \in \mathbb{R}$ 
denotes the forecast value of the ensemble member ${\bf z}_i^{\rm f} \in \mathbb{R}^{40}$ 
at grid point $x_k$ and ${\bf D}(x_k) = \{d_{ij}(x_k)\} \in \mathbb{R}^{M\times M}$.
 
The results for the hybrid NETF-LETKF filter and the hybrid second-order corrected 
ETPF-LETKF can be found in Figure \ref{figure4.1}. The hybrid second-order corrected 
ETPF-LETKF is implemented using the Sinkhorn approximation with $\lambda = 10$ and 
leads to significant improvements over the hybrid NETF-LETKF and also over the hybrid 
ETPF-LETKF of \cite{sr:CRR15}. The CRPS leads to a qualitatively similar assessment. 


\subsection{Estimating parameters for a dynamic scene viewing model} \label{sec:num3}

The scene-viewing model {\it SceneWalk}, as recently proposed by \cite{sr:Engbert2015},
provides a relatively simple mathematical model for a sequence of eye fixations during
scene viewing. The model  dynamically evolves a two-dimensional array of probabilities, $\pi_{ij}(t)$, 
for the next fixation target, which is conditioned on past fixations. More specifically,
the model consists of two sets of ordinary differential equations
\begin{eqnarray} \label{SW1}
\frac{{\rm d}A_{ij}(t)}{{\rm d}t} &=& -\omega_A A_{ij}(t) + \omega_A \frac{S_{ij} \cdot G_A(x_i,y_j;x_f,y_f) }{\sum_{kl} S_{kl} \cdot G_A(x_k,y_l;x_f,y_f)} \\
 \frac{{\rm d}F_{ij}(t)}{{\rm d}t} &=& -\omega_F F_{ij}(t) + \omega_F \frac{G_F(x_i,y_j;x_f,y_f)}{\sum_{kl} G_F(x_k,y_l;x_f,y_f)}  \; \label{SW2}
\end{eqnarray}
for the spatial attention and fixation, respectively, together with
a set of transformation rules
\begin{equation}
u_{ij}(t) = \frac{[A_{ij}(t)]^\lambda}{\sum_{kl}[A_{kl}(t)]^\lambda}-c_{inhib}\frac{[F_{ij}(t)]^\gamma}{\sum_{kl}[F_{kl}(t)]^\gamma}  \;,
\end{equation}
\begin{equation}
u^*(u) = \left\{ 
\begin{array}{ll}
u  & u>\eta \\
\eta e^{\frac{u-\eta}{\eta}} & u\leq \eta
\end{array}
\right. \;,
\end{equation}
which finally produce the desired array of fixation probabilities
\begin{equation}
\pi_{ij}(t) = (1-\zeta)\frac{u^*_{ij}(t)}{\sum_{kl} u^*_{kl}(t)} + \zeta \frac{1}{\sum_{kl} 1} \;.
\end{equation}
The functions $G_{A/F}$ in (\ref{SW1})-(\ref{SW2}) are Gaussians given by
\begin{equation}
\label{eq:Gaussians}
G_{A/F}(x,y; x_f, y_f) = \frac{1}{2\pi\sigma_{A/F}^2}\exp\left(-\frac{(x-x_f)^2+(y-y_f)^2}{2 \sigma_{A/F}^2}\right) 
\end{equation}
and $\{S_{ij}\}$ is a static saliency map.
See \cite{sr:Engbert2015,sr:Schuett2016} for a detailed description of the model. The {\it SceneWalk} 
model contains 9 parameters, which have been estimated in \cite{sr:Schuett2016} using
maximum likelihood estimates. Here we estimate $\sigma_F$ in (\ref{eq:Gaussians}) and 
$\omega_F$ in (\ref{SW2}) with the remaining seven parameter values taken from \cite{sr:Schuett2016}. We start from a uniform prior 
over the interval $[1,5]$ for the first variable and a uniform prior over the interval $[8,16]$ for 
the second variable, respectively. 

Our experiments consist of first computing the importance weights for each sample from the
prior under given pool of fixation paths and then using an LETF to transform those samples into
equally weighted samples from the posterior parameter distribution. We wish to demonstrate the
impact of different LETFs in terms statistical consistency and distribution of their posterior
samples. 

\begin{figure}
\begin{center}
\includegraphics[width=0.47\textwidth]{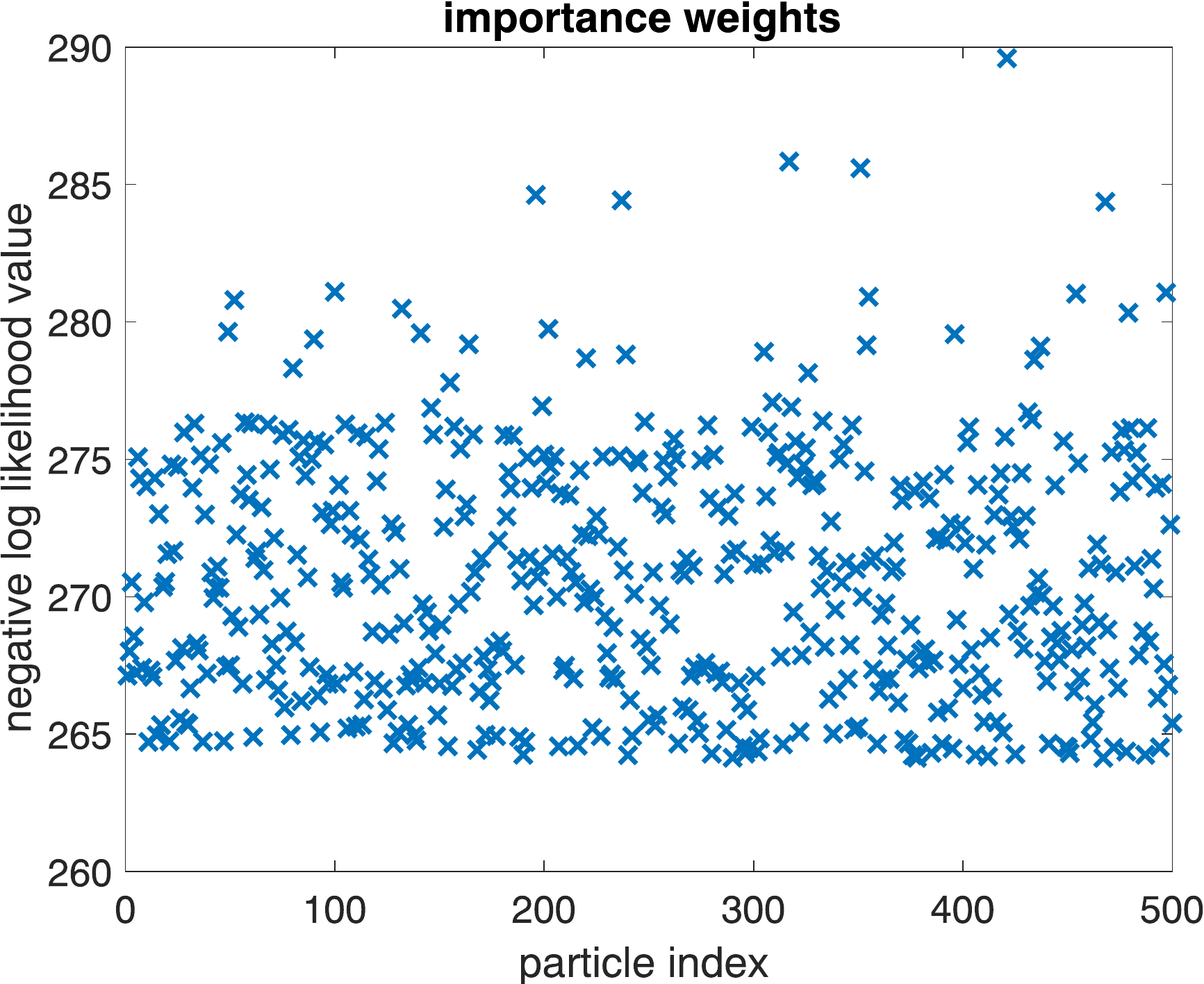}
\end{center}
\caption{Importance weights for $M=500$ samples in two-dimensional parameter space for the {\it Scene Walk} model. The effective
sample size is $M_{\rm eff} \approx 90$.}
\label{figure3.1}
\end{figure}

The importance weights resulting from a given pool of scan paths and $M=500$ samples from
the prior distribution can be found in Figure \ref{figure3.1}. The effective sample size is
about ninety. 

\begin{figure}
\begin{center}
\includegraphics[width=0.47\textwidth]{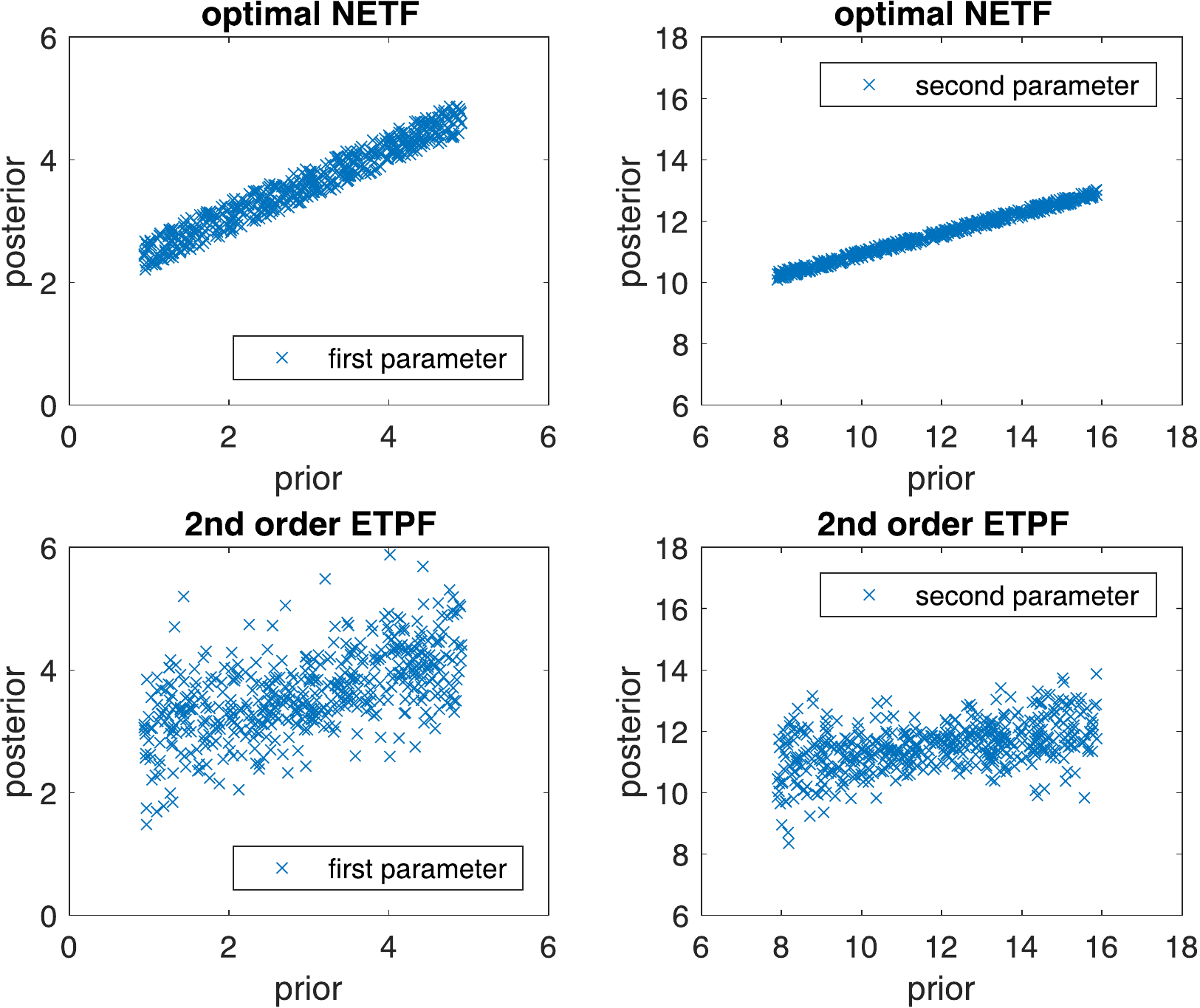}\qquad
\includegraphics[width=0.47\textwidth]{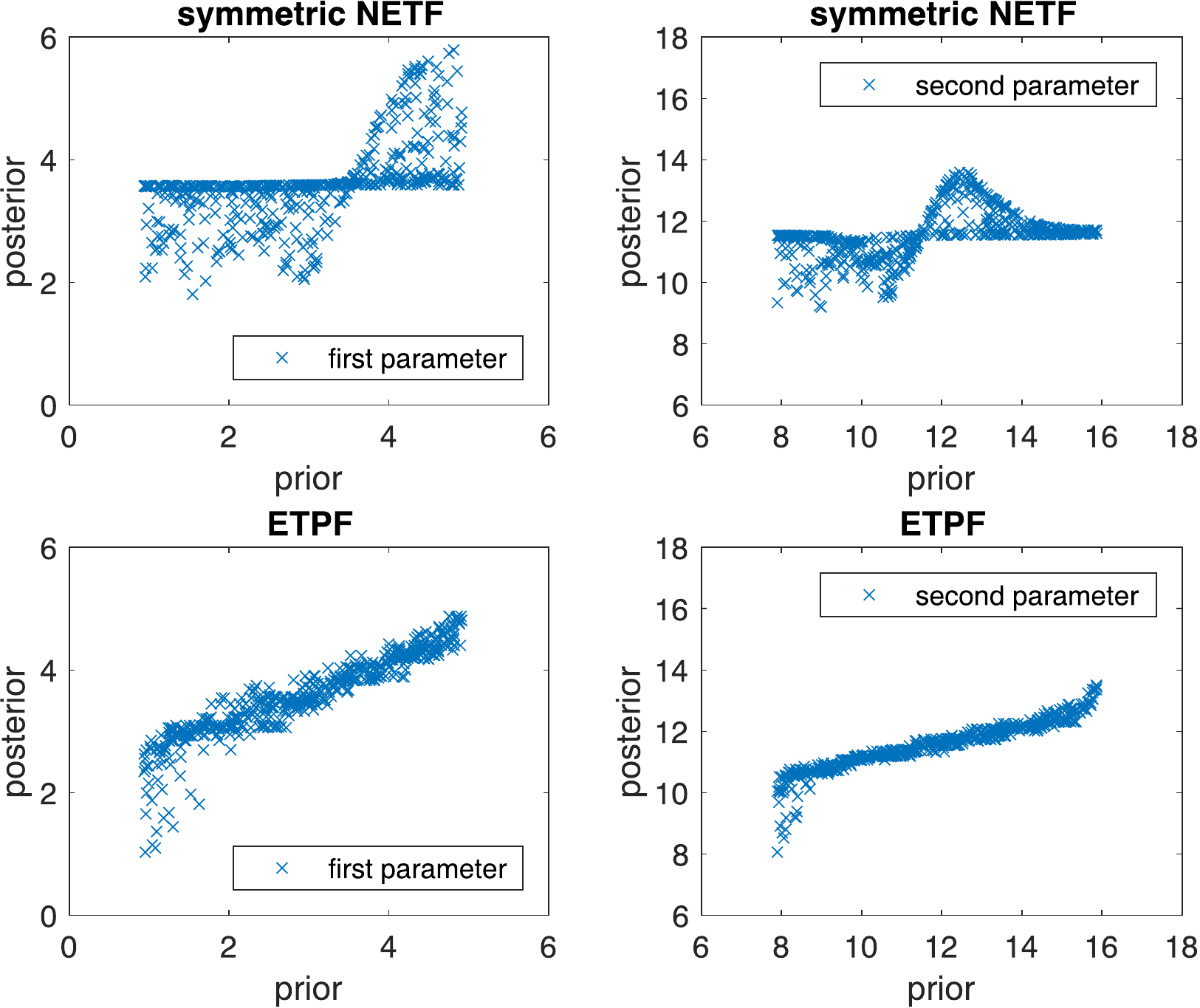}
\end{center}
\caption{Prior vs posterior samples for {\it Scene Walk model}: optimal NEFT (left panel, top row), symmetric NETF (right panel, top row). The two panels also show the ETPF (right panel, bottom row) and the 2nd order 
corrected ETPF (left panel, bottom row) for comparison. Both the optimal NETF and
the ETPF lead to relatively concentrated sample sets, following nearly liner relationships.} \label{figure3.2}
\end{figure}

\begin{figure}
\begin{center}
\includegraphics[width=0.47\textwidth]{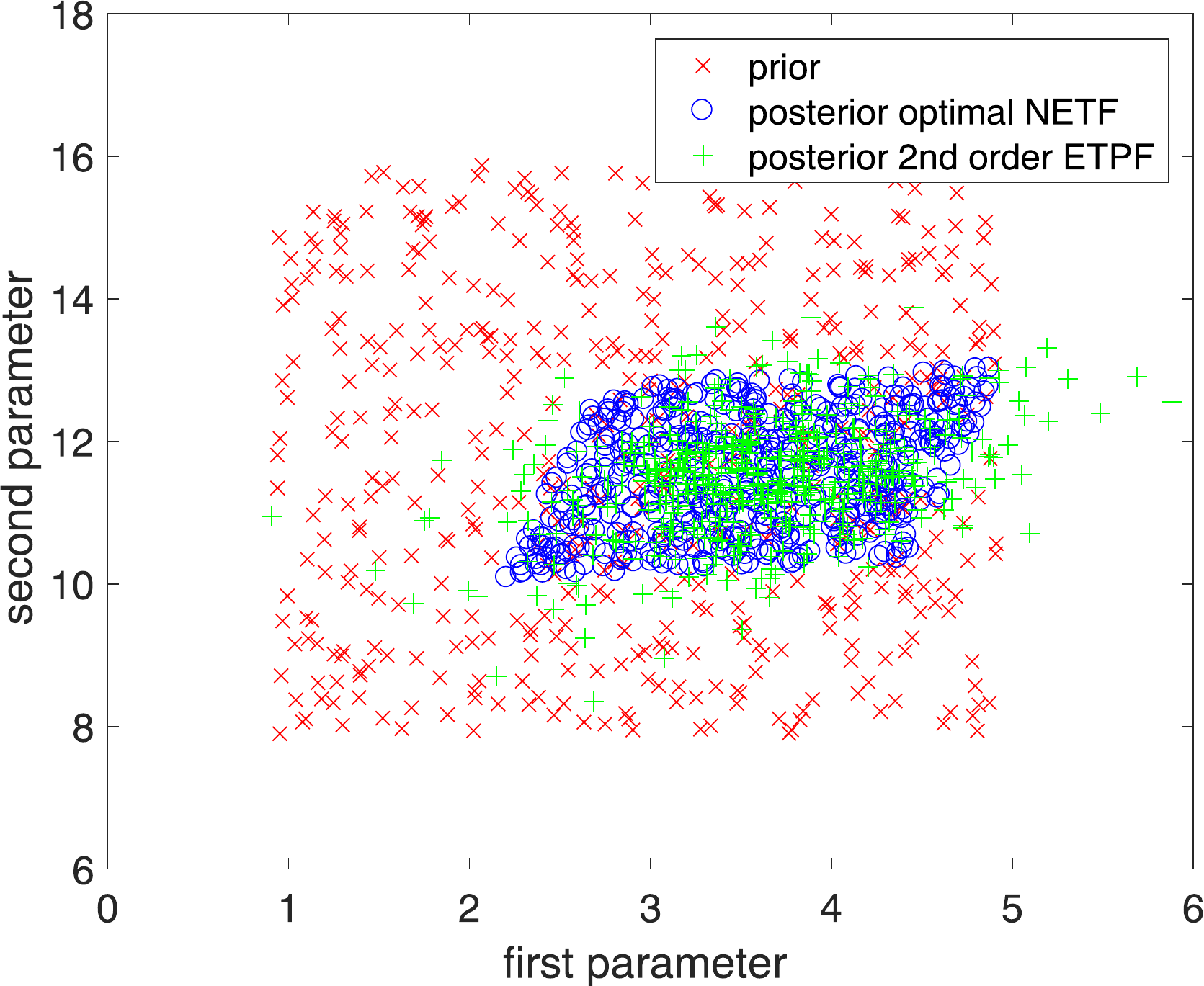}\qquad
\includegraphics[width=0.47\textwidth]{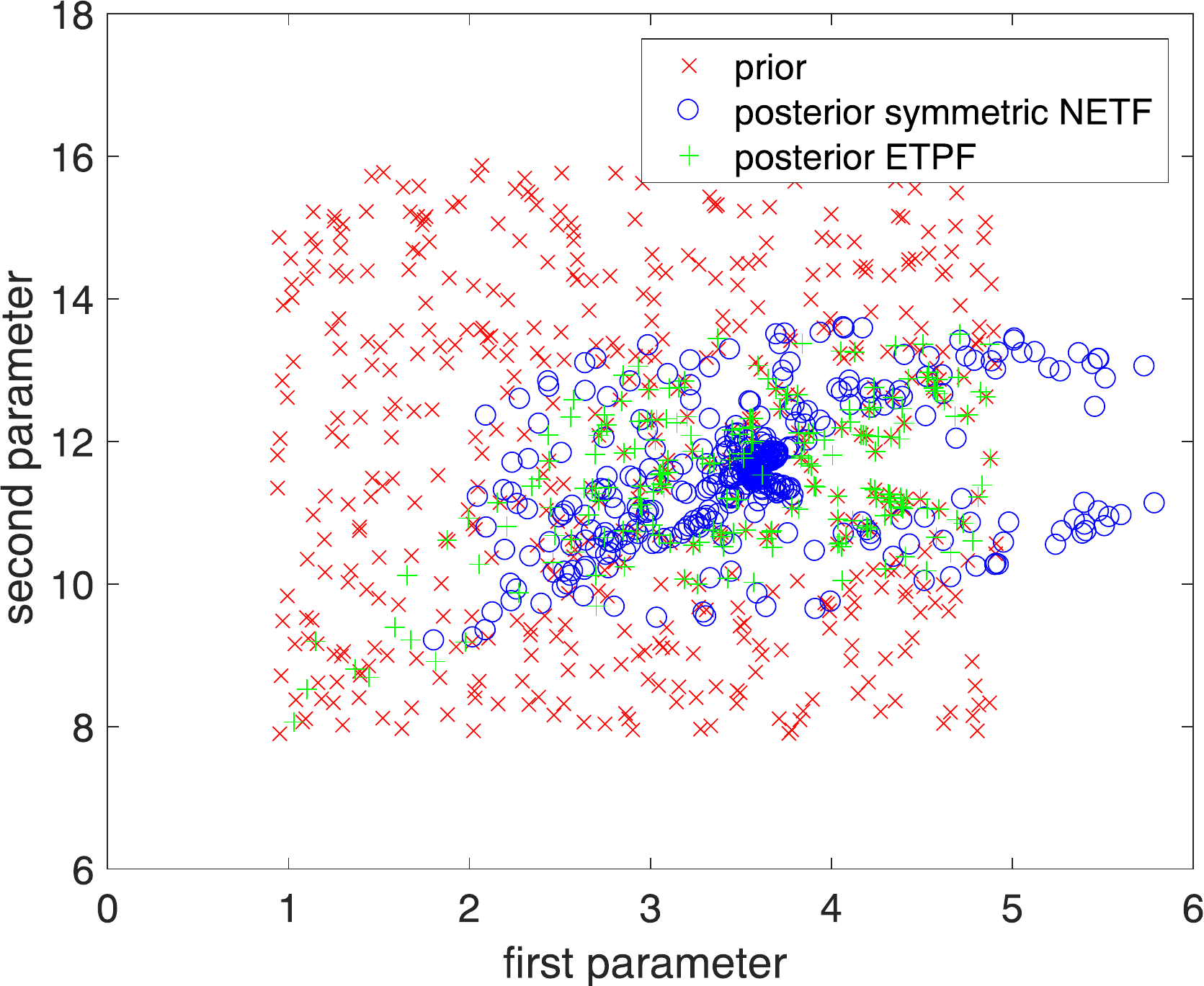}
\end{center}
\caption{Prior vs posterior samples for {\it Scene Walk} model: optimal NEFT (left panel), symmetric NETF (right panel).
Both panels show the ETPF and the 2nd-order corrected ETPF, respectively, 
for comparison. It can be clearly seen that the symmetric 
NETF and the 2nd order corrected ETPF lead to posterior samples which are outside 
the range of the prior samples.} \label{figure3.3}
\end{figure}

We implement the NETF method with ${\bf Q} = {\bf I}$ (symmetric NETF), the NETF with
the optimal ${\bf Q}$ (optimal NETF), the ETPF, and the second-order accurate ETPF.
The distribution of transformed versus prior sample values for each of the two parameters 
separately can be found in Figure \ref{figure3.2}. While the optimal NETF leads to a nearly
linear relation between the prior and transformed samples, the symmetric NETF leads to a rather
non-regular structure. At the same time we find that the second-order accurate ETPF leads to
large fluctuations in the transformed samples with some samples leaving the prior range. 
Since this behavior is violating Bayes' law, it must be seen as a undesirable effect of enforcing 
strict second-order accuracy. The associated two-dimensional scatter plots of the prior and
transformed samples can be found in Figure \ref{figure3.3}. These plots show even more clearly
that second-order accurate methods can lead to transformed samples, which violate Bayes' law.
Nevertheless, all methods consider qualitatively capture the posterior distribution.

\section{Conclusions} 

We have proposed and tested second-order variants of the ETPF. These modifications are computationally attractive since it allows one to replace the computationally expensive solution 
of an optimal transport problem by its Sinkhorn approximation. 
Furthermore, if the regularization parameter, $\lambda$, in the
Sinkhorn approximation is set to zero, then we recover the NETF \cite{sr:toedter15} with an optimally chosen orthogonal matrix ${\bf Q}_{\rm opt}$ in (\ref{optimal}, while $\lambda \to
\infty$ leads formally back to the optimal transport implementation of the ETPF. 
As a byproduct, we also found that
the NETF with an optimally chosen orthogonal matrix, ${\bf Q}$, leads to smaller
RMSEs compared to a random choice, as suggested in \cite{sr:toedter15}.

The second-order accurate ETPF can be 
put into the hybrid ensemble transform particle framework of \cite{sr:CRR15} and can be combined
with localization as necessary for spatially extended evolution equation \cite{sr:evensen,sr:reichcotter15,sr:reich15} such as the Lorenz-96 model. 

The numerical findings for the Lorenz-63 and Lorenz-96 models confirm that 
the methodology proposed in this paper together with the hybrid approach of \cite{sr:CRR15} 
provides a powerful framework for performing sequential data assimilation. We mention that
all methods considered in this paper can be combined with alternative proposal densities, 
which lead to more balanced importance weights (\ref{nweights}) \cite{sr:leeuwen15}. 

It should be noted though, that second-order accuracy comes at a price, i.e., the entries of the transformation
matrix $\widehat{\bf D}$ are not necessarily non-negative, as it is the case for the ETPF transformation matrix ${\bf D}$. 
Hence the analysis ensemble is not necessarily contained in the convex hull spanned by the forecast ensemble. This can cause
non-physical states if, for example, the states should only take values in a bounded interval or semi-interval, as has been demonstrated for the  {\it SceneWalk} model.


\section*{Acknowledgments} We like to thank Hans-Rudolf K\"unsch and Sylvain Robert
for discussions on second-order corrections to linear ensemble transform filters. We also thank Ralf Engbert and Heiko Sch\"utt for providing
the data set used in Section \ref{sec:num3}. This research has been partially funded by Deutsche Forschungsgemeinschaft (DFG) through grant CRC 1114 
"Scaling Cascades in Complex Systems", Project (A02) "Multiscale data and asymptotic model assimilation for atmospheric flows".


\bibliographystyle{siam}
\bibliography{survey}

\end{document}